\documentclass{article}
\usepackage{wrapfig}
\usepackage{calc}
\usepackage{amsmath}
\usepackage{amssymb}
\usepackage{fullpage}
\usepackage{epsfig}

\newcommand{\noun}[1]{\textsc{#1}}
\providecommand{\tabularnewline}{\\}

\begin{document}

\title{Accelerated Optimization in the PDE Framework:\\
 Formulations for the Active Contour Case}

\author{Anthony Yezzi\thanks{School of Electrical and Computer Engineering, Georgia Institute of
Technology} and Ganesh Sundarmoorthi\thanks{Electrical Engineering, King Abdullah University of Science and Technology }}

\date{}

\maketitle
Following the seminal work of Nesterov, accelerated optimization methods
have been used to powerfully boost the performance of first-order,
gradient-based parameter estimation in scenarios where second-order
optimization strategies are either inapplicable or impractical. Not
only does accelerated gradient descent converge considerably faster
than traditional gradient descent, but it also performs a more robust
local search of the parameter space by initially overshooting and
then oscillating back as it settles into a final configuration, thereby
selecting only local minimizers with a basis of attraction large enough
to contain the initial overshoot. This behavior has made accelerated
and stochastic gradient search methods particularly popular within
the machine learning community. In their recent PNAS 2016 paper, Wibisono,
Wilson, and Jordan demonstrate how a broad class of accelerated schemes
can be cast in a variational framework formulated around the Bregman
divergence, leading to continuum limit ODE's. We show how their formulation
may be further extended to infinite dimension manifolds (starting
here with the geometric space of curves and surfaces) by substituting
the Bregman divergence with inner products on the tangent space and
explicitly introducing a distributed mass model which evolves in conjunction
with the object of interest during the optimization process. The co-evolving
mass model, which is introduced purely for the sake of endowing the
optimization with helpful dynamics, also links the resulting class
of accelerated PDE based optimization schemes to fluid dynamical formulations
of optimal mass transport.

\section{Introduction}

Following the seminal work of Nesterov, accelerated optimization methods
(sometimes referred to as momentum methods) have been used to powerfully
boost the performance of first-order, gradient-based parameter estimation
in scenarios where second-order optimization strategies are either
inapplicable or impractical. Not only does accelerated gradient descent
converge considerably faster than traditional gradient descent, but
it also performs a more robust local search of the parameter space
by initially overshooting and then oscillating back as it settles
into a final configuration, thereby selecting only local minimizers
with a basis of attraction large enough to contain the initial overshoot.
This behavior has made accelerated and stochastic gradient search
methods particularly popular within the machine learning community
\cite{Mukherjee13,Li15,Krichene15,Jojic10,Ji09,Hu09,Ghadimi16,Flammarion15,Bubeck15,ODonoghue15}.
So far, however, accelerated optimization methods have been restricted
to searches over finite dimensional parameter spaces. 

Recently, however, Wibisono, Wilson, and Jordan outlined a variational
ODE framework in \cite{Wibisono16} (which we will summarize briefly
in Section \ref{sec:ode}) formulated around the Bregman divergence
and which yields the continuum limit of a broad class of accelerated
optimization schemes, including that of Nesterov's accelerated gradient
method \cite{Nesterov83} whose continuum ODE limit was also demonstrated
by Su, Boyd, and Candes in \cite{Su14}. Here he will show how a similar
high level framework may be adapted for infinite dimensional manifolds
through the formulation of a generalized time-explicit action which
can be viewed both as a specialization and generalization\footnote{Since we abandon the more general Bregman divergence in favor of simpler
inner products, which, however, depend on the more general structure
of the tangent space for the associated infinite dimensional manifold.} of the Bregman Lagrangian presented in \cite{Wibisono16}. While
the extension we outline from the ODE framework into the PDE framework
is general enough to be applied to a variety of infinite-dimensional
or distributed-parameter optimization problems (dense shape reconstruction/inversion,
optical flow estimation, image restoration, etc.) the specific examples
presented here will focus on the active contour and active surface
based optimization. 

Moving into the infinite dimensional framework introduces additional
mathematical, numerical, and computational challenges and technicalities
which do not arise in finite dimensions. For example, the evolving
parameter vector in finite dimensional optimization can naturally
be interpreted as a single moving particle in $\mathbf{R}^{n}$ with
a constant mass which, in accelerated optimization schemes, gains
momentum during its evolution. Since the mass is constant and fixed
to a single particle, there is no need to explicitly model it. When
evolving a continuous curve, surface, region, or function, however,
the notion of accumulated momentum during the acceleration process
is much more flexible, as the corresponding conceptual mass can be
locally distributed in several different ways throughout the domain
which will in turn significantly affect the evolution dynamics. In
fact we intend to exploit this added design flexibility to further
capture some of the same coarse-to-fine regularization properties
of Sobolev gradient flows \cite{Yang15,Sundaramoorthi06} within the
accelerated optimization context as well, but with far less computational
cost.

The discrete implementation of accelerated PDE models will also differ
greatly from existing momentum based gradient descent schemes in finite
dimensions. Spatial and temporal steps sizes will be determined based
on CFL stability conditions for finite difference approximations of
the PDE's. Finally, in the PDE framework, viscosity solutions will
be required in most cases to propagate through shocks and rarefactions
that may occur during the evolution of a continuous front, a phenomenon
which manifests itself differently and is therefore handled differently
in the finite dimensional case. As such, these considerations will
also impact the numerical discretization of accelerated PDE models.

Finally, in part due to these different discretization criteria and
in part to avoid unnecessary complexity in the manifold case, we will
abandon the Bregman Lagrangian described in \cite{Wibisono16} and
will instead exploit a simpler time-explicit\emph{ generalized action}
which will allow us to work directly with the continuum velocity of
the evolving entity rather than finite displacements with the Bregman
divergence. Especially for the case of curves and surfaces considered
here, this avoids the complication of calculating geodesic distances
on highly curved, infinite-dimensional manifolds, but lets us work
more easily in the tangent space instead.

\section{Background and Prior Work}

Geometric partial differential equations have played an important
role in image analysis and computer vision for several decades now.
Applications have ranged from low-level processing operations such
as denoising using anisotropic diffusion, blind deconvolution, and
contrast enhancement; to mid-level processing such as segmentation
using active contours and active surfaces, image registration, and
motion estimation via optical flow; to higher level processing such
as multiview stereo reconstruction, visual tracking, SLAM, and shape
analysis. See, for example, \cite{Sethian-book,Sapiro-book,OsherParagios-book}
for introductions to PDE methods already established within computer
vision within the 1990's, including level set methods \cite{OsherSethian88}
already developed in the 1980's for shape propagation. Several such
PDE methods have been formulated, using the calculus of variations
\cite{Troutman96} as gradient descent based optimization problems
in functional spaces, including geometric spaces of curves and surfaces. 

During the past decade a popular trend has arisen whereby several
such variational problems, which are non-convex, have been reformulated
and relaxed to convex optimization problems \cite{Chan06,Chambolle11,Goldstein10,Pock09},
which allows one to build on the wealth of algorithms developed in
the optimization literature \cite{Boyd04}. While such methods have
led to efficient and robust numerical schemes, the class of problems
for which such reformulations apply are a limited class. We seek to
develop optimization methods for a wider class of (non-convex) problems. 

Recently, Chaudhari \emph{et. al.} have established connections between
relaxation techniques used in training deep neural networks, and PDE's
in\emph{ \cite{Chaudhari17deep}} based on the continuum Fokker-Planck
equation limit. They, in turn, develop and demonstrate improved implementations
of stochastic gradient descent based on the viscous Hamilton-Jacobi
equation. Subsequently, in \cite{Chaudhari17stochastic}, Chaudhari
and Soatto demonstrate that stochastic gradient descent (SGD) methods
perform variational inference (although not on the original loss function).
While they do exploit momentum to accelerate convergence in their
numerical algorithms, this acceleration component is introduced on
the backend of the final discreet algorithm. The methodology presented
in Section \ref{sec:conserved}, through the incorporation of an auxiliary
evolving density function, offers a potential strategy to directly
integrate acceleration into their original continuum PDE formulation
of SDG as well. However, our focus here will remain exclusively on
acceleration, by itself, within the continuum PDE framework. 

\subsection{Geometric Active Contours (an example of gradient PDE optimization)\label{sec:active-contours}}

For example, several active contour models are formulated as gradient
descent PDE flows of application-specific energy functionals $E$
which relate the unknown contour $C$ to given data measurements.
Such energy functionals are chosen to depend only upon the geometric
shape of the contour $C$, not its parameterization. Under these assumptions
the first variation of $E$ will have the following form
\begin{align}
\delta E & =-\int_{C}f\,(\delta C\cdot N)\,ds\label{eq:gradient}
\end{align}
where $fN$ represents a perturbation field along the unit normal
$N$ at each contour point and $ds$ denotes the arclength measure.
Note that the first variation depends only upon the normal component
of a permissible contour perturbation $\delta C$. The form of $f$
will depend upon the particular choice of the energy. For example,
in the popular Chan-Vese active contour model \cite{Chan01} for image
segmentation, $f$ would be expressed by $(I-c_{1})^{2}-(I-c_{2})^{2}+\alpha\kappa$
where $I$ denotes the image value at a given contour point, $\alpha$
an arclength penalty weight, $\kappa$ the curvature at a given contour
point, and $c_{1}$ and $c_{2}$ the means of the image inside and
outside the contour respectively. As an alternative example, the geodesic
active contour model \cite{Caselles98,Kich-ARMA} would correspond
to $f=\phi\kappa N-(\nabla\phi\cdot N)N$ where $\phi>0$ represents
a point measurement designed to be small near a boundary of interest
and large otherwise. In all cases, though, the gradient descent PDE
will the following explicit form.
\begin{equation}
\frac{\partial C}{\partial t}=fN\quad\mbox{[explicit gradient flow]}\label{eq:gradient-flow}
\end{equation}
This class of contour flows, evolving purely in the normal direction,
may be implemented implicitly in the level set framework \cite{OsherSethian88}
by evolving a function $\psi$ whose zero level set represents the
curve $C$ as follows
\[
\frac{\partial\psi}{\partial t}=-\hat{f}\|\nabla\psi\|\quad\mbox{[implicit level set flow]}
\]
where $\hat{f}(x,t)$ denotes a spatial extension of $f(s,t)$ to
points away from the curve.

\subsection{Sobolev gradients for more robust coarse-to-fine PDE based optimization\label{sec:sobolev}}

The most notorious problem with most active contour and active surface
models is that the normal speed function $f$ depends pointwise upon
noisy or irregular data measurements, causing immediate fine scale
perturbations in the evolving contour which cause it to become very
easily attracted to (and trapped within) spurious local minimizers.
This often makes the active contour model strongly dependent upon
initialization, except for a limited class of convex or poly-convex
energy functionals for which numerical schemes can be devised to reach
global minimzers reliably. The traditional way to combat this sensitivity
is to add strong regularizing terms to the energy functional which
penalize fine scale irregularities in the contour shape. Similar problems
and regularization strategies are applied in other PDE based optimization
applications outside the realm of the illustrative active contour
example being considered here (for example, in Horn and Schunck style
optical flow computation \cite{Horn81}).

This energy regularization strategy has two drawbacks. First, most
regularizers lead to second order (or higher) diffusion terms in the
gradient contour flow, which impose much smaller time step limitations
on the numerical discretization of the evolution PDE. Thus, significantly
more evolution steps are required, which incurs a heavy computational
cost in the minimization process. Second, regularizers, while endowing
a level of resistance to noise and spurious structure, impose regularity
on the final converged contour as well, making it difficult or impossible
to capture features such as sharp corners, or narrow protrusions/inlets
in the detected shape. This can lead to unpleasant trade-offs in several
applications.

For the illustrative case of active contours, significantly improved
robustness in the gradient flow, without additional energy regularization,
can be attained by using geometric Sobolev gradients \cite{Charpiat2005,Charpiat2005_Shape_metrics}\cite{Sundaramoorthi05,Sundaramoorthi07}
in place of the standard $L^{2}$-style gradient used in traditional
active contours. We refer to this class of active contours as Sobolev
active contours, whose evolution may be described by the following
integral-partial-differential equation
\begin{equation}
\frac{\partial C}{\partial t}=\left(fN\right)\ast K\quad\mbox{[Sobolev gradient flow]}\label{eq:sobolev-flow}
\end{equation}
Here $\ast$ denotes convolution in the arclength measure with a smoothing
kernel $K$ to invert the linear Sobolev gradient operator. The numerical
implementation is not carried out this way, but the expression gives
helpful insight into how the Sobolev gradient flow (\ref{eq:sobolev-flow})
relates to the usual gradient flow (\ref{eq:gradient-flow}). Namely,
the optimization process (rather than the energy functional itself)
is regularized by averaging point-wise gradient forces $fN$ through
the kernel $K$ to yield a smoother contour evolution. This does not
change the local minimizers of the energy functional, nor does it
impose extra regularity at convergence, but it induces a coarse-to-fine
evolution behavior \cite{Sundaramoorthi06,Yang15} in the contour
evolution, making it much more resistant to spurious local minima
due to noise or other fluctuations in $f$.

\begin{wrapfigure}[21]{o}{0.62\columnwidth}%
\vspace{-8mm}
\hspace{-4mm}%
\begin{tabular}{ccccc}
\multicolumn{5}{c}{{\small{}Active contour (evolving left-to-right) without regularization}}\tabularnewline
\epsfig{figure=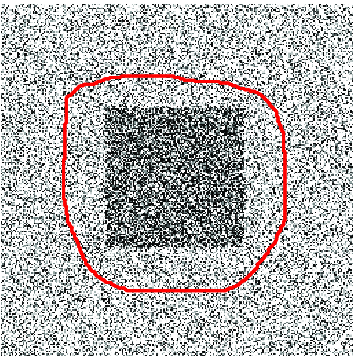,width=.9in} \hspace{-7mm} & \epsfig{figure=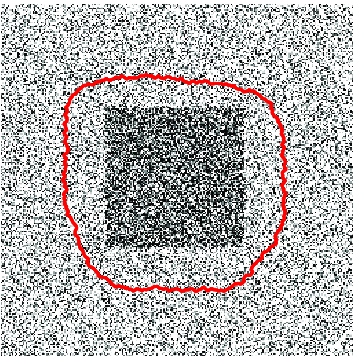,width=.9in}\hspace{-7mm}  & \epsfig{figure=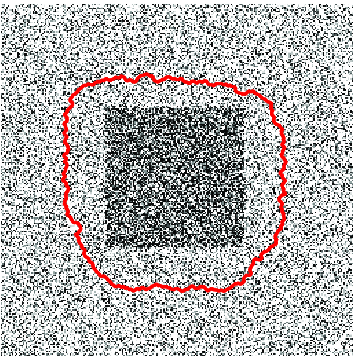,width=.9in} \hspace{-7mm} & \epsfig{figure=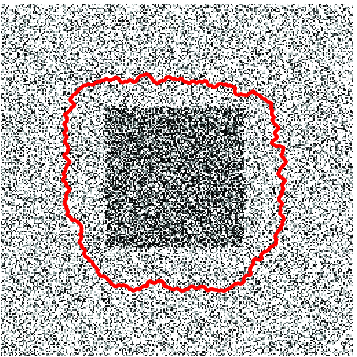,width=.9in} \hspace{-7mm} & \epsfig{figure=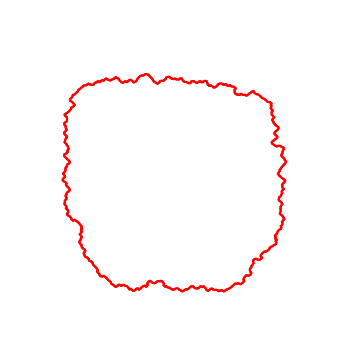,width=.9in} \tabularnewline
\vspace{-4mm}
 &  &  &  & \tabularnewline
\multicolumn{5}{c}{{\small{}Active contour (evolving left-to-right) with added regularization}}\tabularnewline
\epsfig{figure=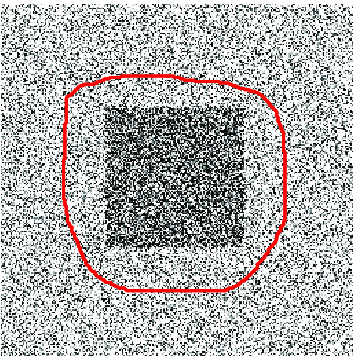,width=.9in}\hspace{-7mm} & \epsfig{figure=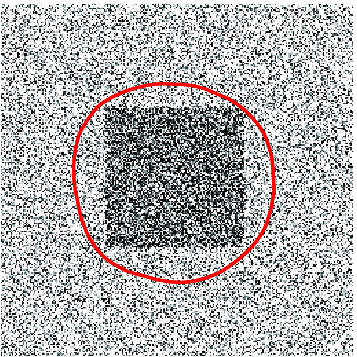,width=.9in} \hspace{-7mm} & \epsfig{figure=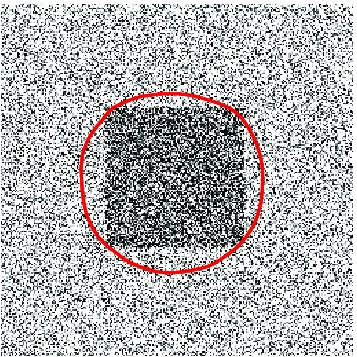,width=.9in} \hspace{-7mm} & \epsfig{figure=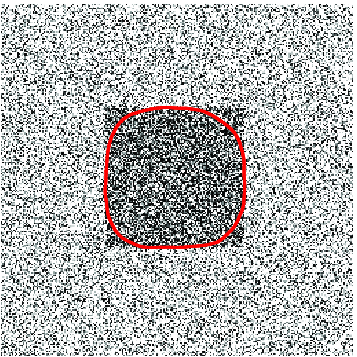,width=.9in} \hspace{-7mm} & \epsfig{figure=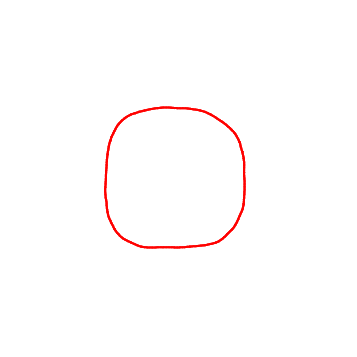,width=.9in} \tabularnewline
\vspace{-4mm}
 &  &  &  & \tabularnewline
\multicolumn{5}{c}{{\small{}Sobolev active contour (evolving left-to-right) without regularization}}\tabularnewline
\epsfig{figure=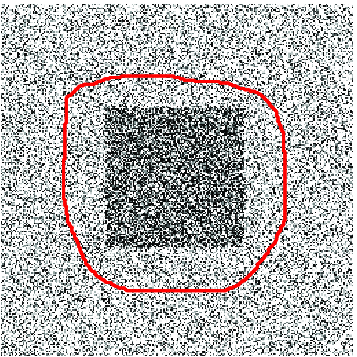,width=.9in} \hspace{-7mm} & \epsfig{figure=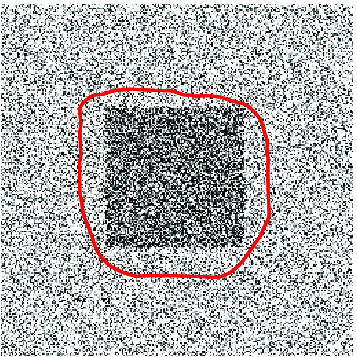,width=.9in} \hspace{-7mm} & \epsfig{figure=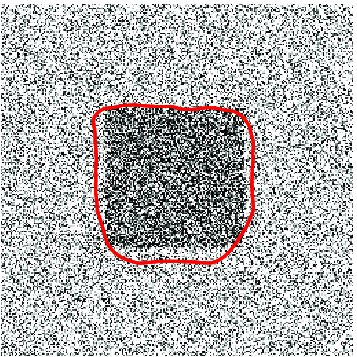,width=.9in} \hspace{-7mm} & \epsfig{figure=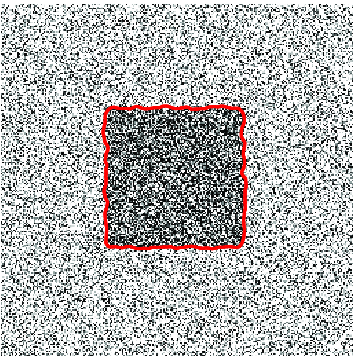,width=.9in} \hspace{-7mm} & \epsfig{figure=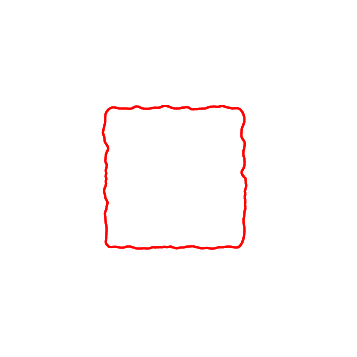,width=.9in} \tabularnewline
\end{tabular}\vspace{-3mm}
\caption{\label{fig:square} Sobolev gradients versus energy regularization}

\end{wrapfigure}%

The regularity of the coarse-to-fine Sobolev gradient flow compared
with regularity imposed on the energy functional is illustrated in
Figure \ref{fig:square}. Along the top row we see the evolution of
a standard active contour in a very noisy image without regularizing
terms in the energy function to keep the contour smooth. The contour
quickly gets trapped in a noisy local minimum configuration before
reaching the desired square boundary. Of course, we can add a regularizing
term to the energy to prefer smoother contours. We see in the middle
row (b) that this fixes the noise problem but does not allow us to
capture the sharp corners of the square. Along the bottom row, instead,
we show the evolution of the Sobolev active contour for the original
{\em unregularized} energy from the top row. The initial stages
of the evolution maintain a smooth contour, not because the Sobolev
gradient prefers a smooth contour, but because it prefers a smooth
evolution. As the Sobolev active contour nears the boundary of the
square, finer scale motions are incorporated to bring out the corners.
The final converged contour responds to local noise, but only in the
vicinity of a desired minimizer.

However, while the Sobolev gradient descent method is extremely successful
in making an active contour or surface (or other evolving classes
of functions) resistant to a large class of unwanted local minimizers,
it comes at heavy computation cost. The spatial integration of gradient
forces along the evolving front must occur during every time step,
and while there are tricks to do this quickly for closed 2D curves
\cite{MennucciNEW,Sundaramoorthi-CDC,Sundaramoorthi-SIAM,Bardelli}
there are no convenient alternatives for 3D surfaces, nor for regions
(even in 2D) when applying Sobolev gradient flows to other functional
objects (images, optical flow, etc.). The linear operator inversion
imposes a notable per-iteration cost, which we will instead distribute
across iterations in the upcoming accelerated coupled PDE evolution
schemes.

\subsection{Momentum methods and Nesterov's Accelerated Gradient}

If we step back to the finite dimensional case, an alternative and
computationally cheaper method to regularize any gradient descent
based iteration scheme is to employ the use of momentum. In such schemes
the new update becomes a weighted combination of the previous update
(the momentum term) and the newly computed gradient at each step.
This leads to a temporal averaging of gradient information computed
and accumulated during the evolution process itself, rather than a
spatial averaging that occurs independently during each time step.
As such it adds insignificant per-iteration computation cost while
significantly boosting the robustness (and often the convergence speed)
of the optimization process.

Momentum methods, including stochastic variants \cite{Ghadimi16,Hu09},
have become very popular in machine learning in recent years \cite{Bubeck15,Flammarion15,Ji09,Jojic10,Krichene15,Mukherjee13,ODonoghue15,Li15}.
Strategic dynamically changing weights on the momentum term can further
boost the descent rate. Nesterov put forth the following famous scheme
\cite{Nesterov83} which attains an optimal rate of order $\frac{1}{t^{2}}$
in the case of a smooth, convex energy function $E(x)$
\[
y_{k+1}=x_{k}-\frac{1}{\beta}\nabla E(x_{k}),\qquad x_{k+1}=(1-\gamma_{k})y_{k+1}+\gamma_{k}y_{k},\qquad\gamma_{k}=\frac{1-\lambda_{k}}{\lambda_{k}+1},\qquad\lambda_{k}=\frac{1+\sqrt{1+4\lambda_{k-1}^{2}}}{2}
\]
where $x_{k}$ is the $k$-th iterate of the algorithm, $y_{k}$ is
an intermediate sequence, and $\gamma_{k}$ are dynamically updated
weights.

\subsection{A Variational Framework for Accelerated ODE Optimization\label{sec:ode}}

Recently in \cite{Wibisono16} Wibisono, Wilson, and Jordan presented
a variational generalization of Nesterov's \cite{Nesterov83} and
other momentum based gradient descent schemes in $\mathbb{R}^{n}$
based on the Bregman divergence of a convex distance generating function
$h$ 
\begin{equation}
D(y,x)=h(y)-h(x)-\left\langle \nabla h(x),y-x\right\rangle \label{eq:breg-divergence}
\end{equation}
and careful discretizations of the Euler-Lagrange equation for the
time integral (evolution time) of the following Bregman Lagrangian
\[
{\cal L}(X,V,t)=e^{a(t)+\gamma(t)}\left[D(X+e^{-a(t)}V,X)-e^{b(t)}{\bf U}(X)\right]
\]
where the potential energy ${\bf U}$ represents the cost to be minimized.
In the Euclidean case, where $D(y,x)=\frac{1}{2}\|y-x\|^{2}$, this
simplifies to 
\[
{\cal L}=e^{\gamma(t)}\left[e^{-a(t)}\underbrace{\frac{1}{2}\|V\|^{2}}_{{\bf T}}-e^{a(t)+b(t)}{\bf U}(X)\right]
\]
where ${\bf T}$ models the kinetic energy of a unit mass particle
in $\mathbb{R}^{n}$. Nesterov's methods \cite{Nesterov83,Nesterov14,Nesterov13,Nesterov08,Nesterov06,Nesterov05}
belong to a subfamily of Bregman Lagrangians with the following choice
of parameters (indexed by $k>0$) 
\[
a=\log k-\log t,\qquad b=k\,\log t+\log\lambda,\qquad\gamma=k\,\log t
\]
which, in the Euclidean case, yields a time-explicit \emph{generalized
action} (compared to the time-implicit standard action ${\bf T}-{\bf U}$
from classical mechanics \cite{Goldstein02}) as follows.
\begin{equation}
{\cal L}=\frac{t^{k+1}}{k}\left({\bf T}-\lambda k^{2}t^{k-2}{\bf U}\right)\label{eq:time-action}
\end{equation}
In the case of $k=2$, for example, the Euler-Lagrange equations for
the integral of this time-explicit action yield the continuum limit
of Nesterov's accelerated mirror descent\cite{Nesterov05} derived
in both \cite{Su14,Krichene15}. 

\section{Accelerated Optimization in the PDE Framework}

We now develop a general strategy, based on a generalization of the
Euclidean case of Wibisono, Wilson, and Jordan's formulation \cite{Wibisono16}
reviewed in Section \ref{sec:ode}, for extending accelerated optimization
into the PDE framework. While our approach will be motivated by the
variational ODE framework formulated around the Bregman divergence
in \cite{Wibisono16}, we will have to address several mathematical,
numerical, and computational considerations which do not need to be
addressed in finite dimenions.

For example, the evolving parameter vector in finite dimensional optimization
can naturally be interpreted as a single moving particle in $\mathbf{R}^{n}$
with a constant mass which, in accelerated optimization schemes, gains
momentum during its evolution. Since the mass is constant and fixed
to a single particle, there is no need to explicitly model it. When
evolving a continuous curve, surface, region, or function, however,
the notion of accumulated momentum during the acceleration process
is much more flexible, as the corresponding conceptual mass can be
locally distributed in several different ways throughout the domain
which will in turn significantly affect the evolution dynamics. We
outline two different mass models in Sections \ref{sec:constant}
and \ref{sec:conserved} as starting points and show how additional
control of the optimization dynamics can be introduced in conjunction
with the more flexible second mass model by considering independent
mass-related potential energy terms in Section \ref{sec:mass-potential}.
In all cases, the outcome of these formulations will be a coupled
system of first-order PDE's which govern the simultaneous evolution
of the continuous unknown (curves in the case considered here), its
velocity, as well as the supplementary density function which describes
the evolving mass.

In addition, as pointed out from the onset, the numerical discretization
of accelerated PDE models will also differ greatly from existing momentum
based gradient descent schemes in finite dimensions. Spatial and temporal
steps sizes will be determined based on CFL stability conditions for
finite difference approximations of the PDE's and viscosity solution
schemes will be required in most cases to propagate through shocks
and rarefactions that may occur during the evolution of a continuous
front. This is part of the reason we replace the more general Bregman-Lagrangian
in \cite{Wibisono16} with the simpler time-explicit\emph{ generalized
action} (\ref{eq:time-action}), together with the additional benefit
that such a choice allows us to work directly with the continuum velocity
of the evolving entity (or other generalizations that are easily defined
within the tangent space of its relevant manifold) rather than finite
displacements utilized by the Bregman divergence (\ref{eq:breg-divergence}). 

\subsection{General Approach\label{sec:general}}

Just as in \cite{Wibisono16}, the energy functional $E$ to be optimized
over the continuous infinite dimensional unknown (whether it be a
function, a curve, a surface, or a diffeomorphic mapping) will represent
the potential energy term ${\bf U}$ in the time-explicit \emph{generalized
action} (\ref{eq:time-action}). Next, a customized kinetic energy
term ${\bf T}$ will be formulated to incorporate the dynamics of
the evolving estimate during the minimization process. Note that just
as the evolution time $t$ would represent an artificial time parameter
for a continuous gradient descent process, the kinetic energy term
will be linked to artificial dynamics incorporated into the accelerated
optimization process. As such, the accelerated optimization dynamics
can be designed completely independently of any potential physical
dynamics in cases where the unknown might be connected with the motion
of real objects. Several different strategies can be explored, depending
upon the geometry of the specific optimization problem, for defining
kinetic energy terms, including various approaches for attributing
artificial mass (both its distribution and its flow) to the actual
unknown of interest in order to boost the robustness and speed of
the optimization process.

\begin{wrapfigure}[17]{O}{0.45\textwidth}%
\vspace{-10mm}
\centerline{\epsfig{figure=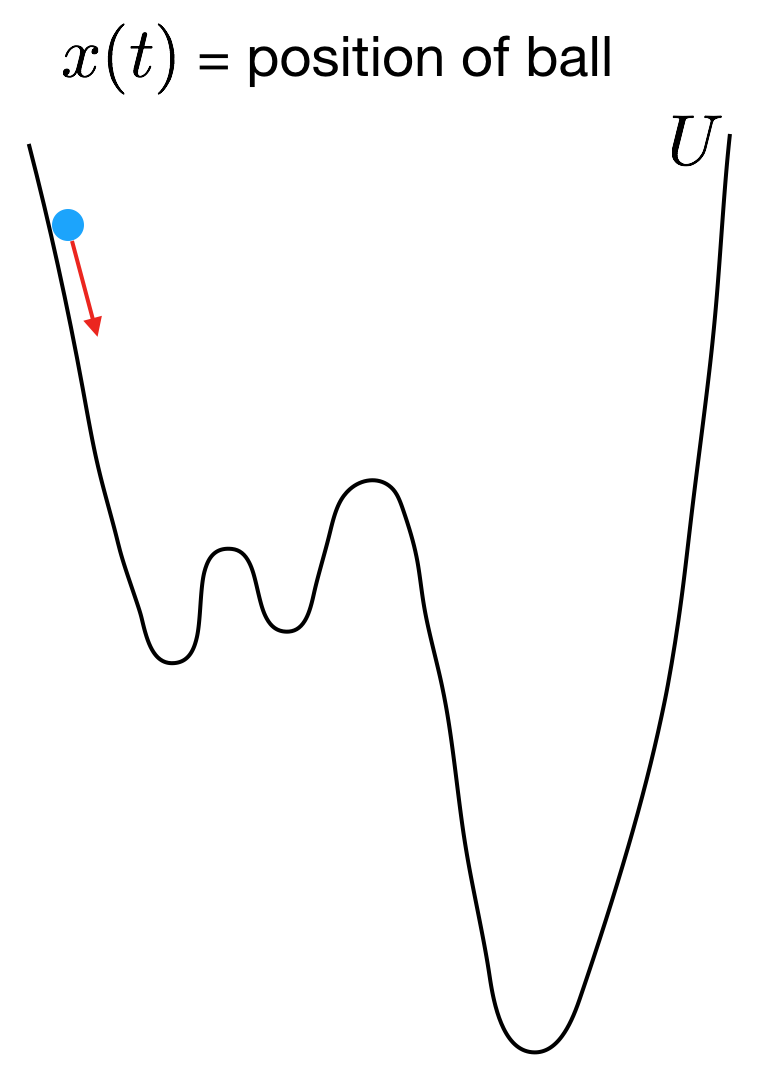,width=2in}} \vspace{-7mm}

{\footnotesize{}\caption{{\small{}Accelerated descent physics interpretation\label{fig:acceleration} }}
}{\footnotesize \par}\end{wrapfigure}%

Once the kinetic energy term has been formulated, the accelerated
evolution will obtained (prior to discretization) using the Calculus
of Variations\cite{Troutman96} as the Euler-Lagrange equation of
the following time-explicit \emph{generalized action integral} 
\begin{equation}
\int\frac{t^{k+1}}{k}\left({\bf T}-\lambda k^{2}t^{k-2}{\bf U}\right)\,dt\label{eq:action-integral}
\end{equation}
In the simple $k=2$ case, the main difference between the resulting
evolution equations versus the classical Principle of Least Action
equations of motion (without the time explicit terms in the Lagrangian)
is an additional friction-style term whose coefficient of friction
decreases inversely proportional to time. This additional term, however,
is crucial to the accelerated minimization scheme. Without such a
frictional term, the Hamiltonian of the system (the total energy ${\bf T}+{\bf U}$),
would be conserved, and the associated dynamical evolution would never
converge to a stationary point. Friction guarantees a monotonic dissipation
of energy, allowing the evolution to converge to a state of zero kinetic
energy and locally minimal potential energy (the optimization objective).

This yields a natural physical interpretation of accelerated gradient
optimization in terms of a mass rolling down a potentially complicated
terrain by the pull of gravity (Figure \ref{fig:acceleration}). In
gradient descent, its mass is irrelevant, and the ball always rolls
downward by gravity (the gradient). As such the gradient directly
regulates its velocity. In the accelerated case, gravity regulates
its acceleration. Friction can be used to interpolate these behaviors,
with gradient descent representing the infinite frictional limit as
pointed out in \cite{Wibisono16}.

Acceleration comes with two advantages. First, whenever the gradient
is very shallow (the energy functional is nearly flat), acceleration
allows the ball to accumulate velocity as it moves so long as the
gradient direction is self reinforcing. As such, the ball approaches
a minimum more quickly. Second, the velocity cannot abruptly change
near a shallow minimum as in gradient descent. Its mass gives it momentum,
and even if the acceleration direction switches in the vicinity of
a shallow minimum, the accumulated momentum still moves it forward
for a certain amount of time, allowing the optimization process to
\emph{look ahead} for a potentially deeper minimizer. 

\subsection{Accelerated Active Contours\label{sec:contours}}

We now illustrate the steps in the process for developing PDE based
accelerated optimization schemes for the specific case of geometric
active contours. The resulting coupled PDE evolutions will retain
the parameterization independent property of gradient descent based
active contours models and will therefore remain amenable to implicit
implementation using Level Set Methods \cite{OsherSethian88}.

We begin, however, by reviewing some basic differential contour evolution
properties that will be useful in deriving accelerated active contour
formulations. In particular, it is useful to understand any contour
evolution behavior in terms of its local geometric frame, consisting
of the unit tangent and normal vectors. 

Let $C(p,t)$ denote an evolving curve where $t$ represents the evolution
parameter and $p\in[0,1]$ denotes an independent parameter along
each fixed curve. The unit tangent, unit normal, and curvature will
be denoted by $T=\frac{\partial C}{\partial s}$, $N$, and $\kappa$
respectively, with the sign convention for $\kappa$ and the direction
convention for $N$ chosen to respect the planar Frenet equations
$\frac{\partial T}{\partial s}=\kappa N$ and $\frac{\partial N}{\partial s}=-\kappa T$,
where $s$ denotes the time-dependent arclength parameter whose derivative
with respect to $p$ yields the parameterization speed $\frac{\partial s}{\partial p}=\left\Vert \frac{\partial C}{\partial p}\right\Vert $.

Letting $\alpha$ and $\beta$ denote the tangential and normal speeds
of the evolving curve\footnote{Note that the instantaneous geometric deformation of the curve is
determined exclusively by the normal speed $\beta$, and that gradient
flows for geometric active contours can all be formulated such that
the tangential speed $\alpha$ vanishes. We will see later that the
same is possible for accelerated flow models as well. },
\begin{equation}
\frac{\partial C}{\partial t}=\alpha T+\beta N\label{eq:velocity}
\end{equation}
the frame itself can be shown to evolve as follows.
\begin{equation}
\frac{\partial T}{\partial t}=\left(\frac{\partial\beta}{\partial s}+\alpha\kappa\right)\,N,\qquad\frac{\partial N}{\partial t}=-\left(\frac{\partial\beta}{\partial s}+\alpha\kappa\right)\,T\label{eq:frame-evolution}
\end{equation}
Differentiating the velocity decomposition (\ref{eq:velocity}) with
respect to $t$, followed by the frame evolution (\ref{eq:frame-evolution})
substitution, yields the acceleration 
\begin{equation}
\frac{\partial^{2}C}{\partial t^{2}}=\left(\frac{\partial\alpha}{\partial t}-\beta\left(\frac{\partial\beta}{\partial s}+\alpha\kappa\right)\right)\,T+\left(\frac{\partial\beta}{\partial t}+\alpha\left(\frac{\partial\beta}{\partial s}+\alpha\kappa\right)\right)\,N\label{eq:acceleration}
\end{equation}
which may be rewritten as the following two scalar evolution equations
for the tangential and normal speeds, in terms of the tangential and
normal components of the contour acceleration, respectively.
\begin{equation}
\frac{\partial\alpha}{\partial t}=\frac{\partial^{2}C}{\partial t^{2}}\cdot T+\beta\left(\frac{\partial\beta}{\partial s}+\alpha\kappa\right),\qquad\frac{\partial\beta}{\partial t}=\frac{\partial^{2}C}{\partial t^{2}}\cdot N-\alpha\left(\frac{\partial\beta}{\partial s}+\alpha\kappa\right)\label{eq:speed-evolution}
\end{equation}

\subsubsection{Contour potential energy\label{sec:potential}}

For geometric active contours, we start by defining the potential
energy ${\bf U}$ to be an originally provided energy functional $E$
which depends only upon the geometric shape of the contour $C$ (not
its parameterization). Under these assumptions the first variation
of the potential energy will have the following form, just as in (\ref{eq:gradient})
presented earlier in Section \ref{sec:active-contours}, where $fN$
denotes the backward local gradient force at each contour point.
\begin{align*}
\delta{\bf U} & =-\int_{C}f\,(\delta C\cdot N)\,ds
\end{align*}

\subsubsection{Constant density model\label{sec:constant}}

To formulate an accelerated evolution model, we define a kinetic energy,
which requires a notion of mass coupled with velocity. The simplest
starting model would be one of constant mass density $\rho$ (per
unit arclength along the contour) and an integral of the squared norm
of the point-wise contour evolution velocity\footnote{A similar kinetic energy model in the context of the classical action
${\bf T}-{\bf U}$, for example, was used to develop dynamic geodesic
snake models for visual tracking in \cite{Niethammer06} }.
\begin{equation}
{\bf T}=\frac{1}{2}\rho\int_{C}\left(\frac{\partial C}{\partial t}\cdot\frac{\partial C}{\partial t}\right)ds\label{eq:kinetic-constant}
\end{equation}
Plugging this into the generalized action integral (\ref{eq:action-integral})
and computing the Euler-Lagrange equation leads to our first, and
simplest, accelerated model.
\begin{equation}
\underbrace{\frac{\partial^{2}C}{\partial t^{2}}}_{\mbox{acceleration}}=\frac{\lambda k^{2}t^{k-2}}{\rho}\underbrace{fN}_{-\mbox{gradient}}-\underbrace{\left(\frac{\partial^{2}C}{\partial s\partial t}\cdot\frac{\partial C}{\partial s}\right)\frac{\partial C}{\partial t}-\frac{\partial}{\partial s}\left(\frac{1}{2}\left\Vert \frac{\partial C}{\partial t}\right\Vert ^{2}\frac{\partial C}{\partial s}\right)}_{\mbox{wave propagation terms}}-\underbrace{\frac{k+1}{t}\frac{\partial C}{\partial t}}_{\mbox{friction}}\label{eq:acceleration-constant}
\end{equation}
If we start with zero initial velocity we can decompose this nonlinear
second-order PDE into the following coupled system of nonlinear first
order PDE's 
\begin{equation}
\frac{\partial C}{\partial t}=\beta N,\qquad\frac{\partial\beta}{\partial t}=\frac{\lambda k^{2}t^{k-2}}{\rho}f+\frac{1}{2}\beta^{2}\kappa-\frac{k+1}{t}\beta\label{eq:accelerated-system-constant}
\end{equation}
Since the contour evolution remains purely geometric (only in the
normal direction $N$) we may also write down an implicit level set
version of the coupled PDE system as follows 
\begin{equation}
\frac{\partial\psi}{\partial t}=\hat{\beta}\|\nabla\psi\|,\qquad\frac{\partial\hat{\beta}}{\partial t}=\frac{\lambda k^{2}t^{(k-2)}}{\rho}\hat{f}+\nabla\cdot\left(\frac{1}{2}\hat{\beta}^{2}\frac{\nabla\psi}{\|\nabla\psi\|}\right)-\frac{k+1}{t}\hat{\beta}\label{eq:levelset-constant}
\end{equation}
where $\hat{f}(x,t)$ and $\hat{\beta}(x,t)$ denote spatial extensions
of $f$ and $\beta$ respectively.

\pagebreak{}

\paragraph{Illustrative results}

\begin{wrapfigure}[35]{O}{0.62\columnwidth}%
\vspace{-4mm}
\hspace{-4mm}%
\begin{tabular}{ccccc}
\epsfig{figure=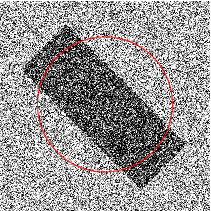,width=.9in} \hspace{-7mm} & \epsfig{figure=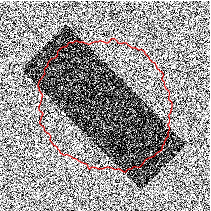,width=.9in}\hspace{-7mm}  & \epsfig{figure=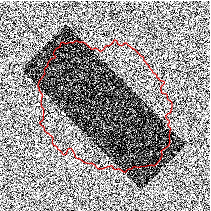,width=.9in} \hspace{-7mm} & \epsfig{figure=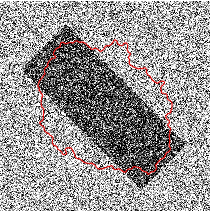,width=.9in} \hspace{-7mm} & \epsfig{figure=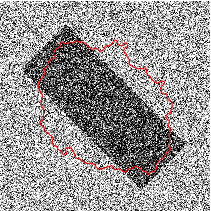,width=.9in} \tabularnewline
\vspace{-5mm}
 &  &  &  & \tabularnewline
\epsfig{figure=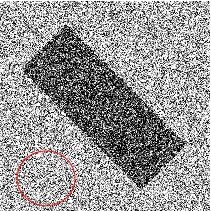,width=.9in}\hspace{-7mm} & \epsfig{figure=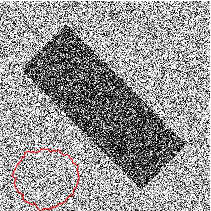,width=.9in} \hspace{-7mm} & \epsfig{figure=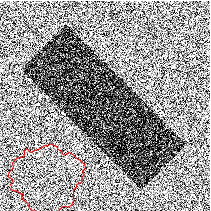,width=.9in} \hspace{-7mm} & \epsfig{figure=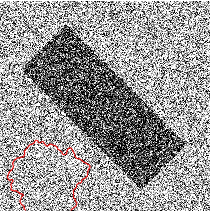,width=.9in} \hspace{-7mm} & \epsfig{figure=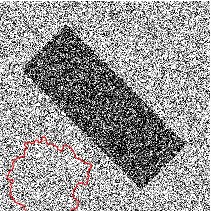,width=.9in} \tabularnewline
\vspace{-5mm}
 &  &  &  & \tabularnewline
\epsfig{figure=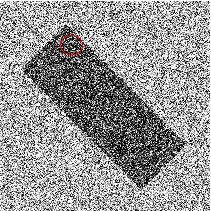,width=.9in} \hspace{-7mm} & \epsfig{figure=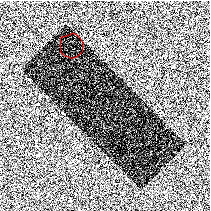,width=.9in} \hspace{-7mm} & \epsfig{figure=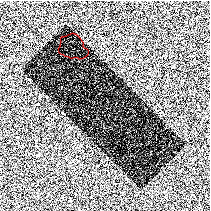,width=.9in} \hspace{-7mm} & \epsfig{figure=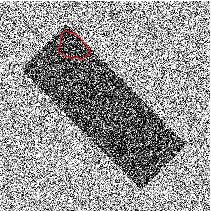,width=.9in} \hspace{-7mm} & \epsfig{figure=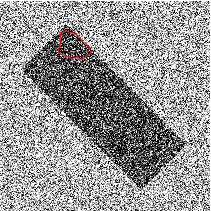,width=.9in}\tabularnewline
\end{tabular}\vspace{-4mm}
\caption{\label{fig:normal}{\small{} Three active contours getting stuck in
different local minima}}
\vspace{1mm}
\hspace{-4mm}%
\begin{tabular}{ccccc}
\epsfig{figure=rect1_init,width=.9in} \hspace{-7mm} & \epsfig{figure=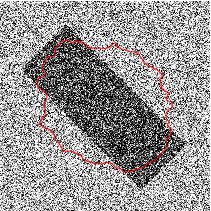,width=.9in}\hspace{-7mm}  & \epsfig{figure=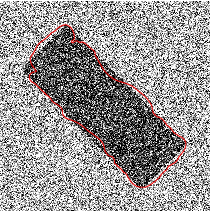,width=.9in} \hspace{-7mm} & \epsfig{figure=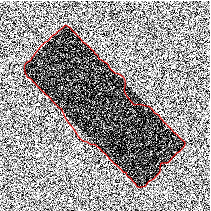,width=.9in} \hspace{-7mm} & \epsfig{figure=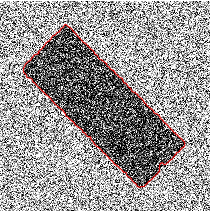,width=.9in} \tabularnewline
\vspace{-5mm}
 &  &  &  & \tabularnewline
\epsfig{figure=rect2_init,width=.9in}\hspace{-7mm} & \epsfig{figure=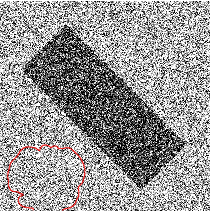,width=.9in} \hspace{-7mm} & \epsfig{figure=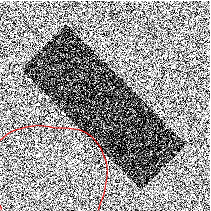,width=.9in} \hspace{-7mm} & \epsfig{figure=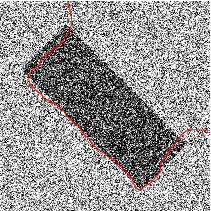,width=.9in} \hspace{-7mm} & \epsfig{figure=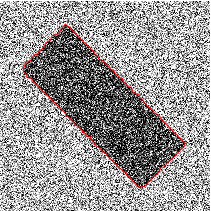,width=.9in} \tabularnewline
\vspace{-5mm}
 &  &  &  & \tabularnewline
\epsfig{figure=rect3_init,width=.9in} \hspace{-7mm} & \epsfig{figure=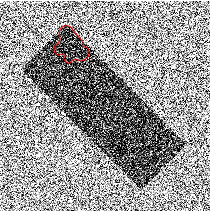,width=.9in} \hspace{-7mm} & \epsfig{figure=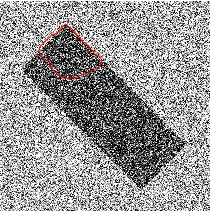,width=.9in} \hspace{-7mm} & \epsfig{figure=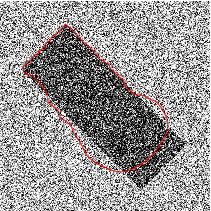,width=.9in} \hspace{-7mm} & \epsfig{figure=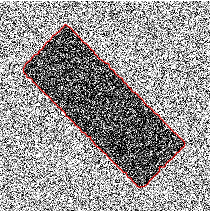,width=.9in}\tabularnewline
\end{tabular}\vspace{-4mm}
\caption{\label{fig:accel}{\small{} Accelerated active contours all converging
to same minimzer}}
\end{wrapfigure}%
The benefits of acceleration can already been seen, even for this
simplest starting kinetic energy model (and for the simplest case
of $k$=2 where the gradient force term remains time-independent)
by comparing Figures \ref{fig:normal} and \ref{fig:accel}. In Figure
\ref{fig:normal} we see three different initial contour placements
(top, middle, bottom) evolving from left-to-right, each getting stuck
in a different local minimizers due to noise, all of which lie very
far away from the desired much deeper minimizer along the rectangle
boundary. Of course, as always, stronger regularizing terms could
be added to the active contour energy functional (or the potential
energy in the accelerated scheme) in order to impose smoothness on
the contour, thereby making it resistant to noise. However, the point
of this synthetic experiment was to create an energy landscape littered
with literally tens of thousands (perhaps even hundreds of thousands)
of local minimizers in order to demonstrate the effects of acceleration.
Furthermore, stronger regularization would come with the additional
sacrifice of being unable to capture the sharp corners of the rectangle
as well as higher computational cost resulting from the smaller step
size constraint in the PDE discretization.

We will not have to make this sacrifice if we regularize the evolution
rather than the contour (as we saw with the Sobolev active contour
in Figure \ref{fig:square}). Indeed, in Figure \ref{fig:accel},
we see the effect of applying this simplest accelerated evolution
scheme (\ref{eq:acceleration-constant}) with the same initial contour
placements and the same energy functional (without additional regularizing
terms). In all three cases, the accelerated PDE system initially pushes
the contour past the noise, driving it toward the deeper minimum along
the rectangle edge. Just as for the Sobolev active contour case in
Figure \ref{fig:square},we see that the accelerated active contour
only captures noisy structure near the desired converged result.

\subsubsection{Conserved flowable mass model \label{sec:conserved}}

The kinetic energy in the accelerated formulation is invented to endow
the minimization evolution process with helpful dynamics for the sake
of faster and more robust convergence. Thus, just as the potential
energy does not actually represent a real physical energy, there is
no need to impose real physical considerations on the kinetic energy
either. Nonetheless, the simple constant density model feels quite
unnatural in that it does not preserve total mass if the contour length
changes during its evolution: mass is created when the contour expands
and is destroyed when the contour contracts. 

A more flexible and natural way to attribute mass to the evolving
contour is to consider an arbitrary and independent distribution of
mass along the contour which evolves as the curve evolves. As such,
the mass density $\rho$ can vary both spatially and temporarily,
while the total integrated mass is still conserved. In such a model,
though, not only does mass evolve as a result of contour shape deformation,
but it may also flow along the contour without changing its geometry
(therefore contributing to the kinetic energy without affecting the
potential energy). A simple interpretation would be that the contour
shape represents a moving container for a fluid which not only gets
pushed around by the extrinsic motion of the container but which may
also flow with an independent relative internal speed $v$ inside
of the container (i.e. along the tangent direction of the contour).
As such, the velocity of each mass particle at a given contour point
would be the sum of the contour velocity and the internal mass flow
velocity.

\[
v=\mbox{internal mass flow speed},\qquad\mbox{total mass velocity}=\frac{\partial C}{\partial t}+v\frac{\partial C}{\partial s}
\]
This suggests a more general kinetic energy model as follows
\begin{equation}
{\bf T}=\int_{C}\frac{1}{2}\rho\left\Vert \frac{\partial C}{\partial t}+v\frac{\partial C}{\partial s}\right\Vert ^{2}ds\label{eq:kinetic-flowable}
\end{equation}
but with the density evolution constrained by the following continuity
equation to ensure local conservation of mass.
\begin{equation}
\underbrace{\frac{\partial\rho}{\partial t}+\frac{\partial}{\partial s}(\rho v)}_{\mbox{mass change}}+\rho\,\underbrace{\left(\frac{\partial^{2}C}{\partial s\partial t}\cdot\frac{\partial C}{\partial s}\right)}_{\mbox{length change}}=0\label{eq:continuity-constraint}
\end{equation}
The latter may be incorporated as a Lagrange multiplier constraint
when computing the Euler-Lagrange equation of the generalized action
integral (\ref{eq:action-integral}). This results in the following
second order PDE's which, together with (\ref{eq:continuity-constraint}),
yield the accelerated system as a coupled evolution of $C$ along
with the auxiliary mass density $\rho$ and internal flow $v$ field
responsible for these helpful dynamics.
\begin{align}
\frac{\partial^{2}C}{\partial t^{2}}\cdot N & =-\left(2v\,\frac{\partial^{2}C}{\partial s\partial t}+v^{2}\frac{\partial^{2}C}{\partial s^{2}}+\frac{k+1}{t}\frac{\partial C}{\partial t}\right)\cdot N+\frac{\lambda k^{2}t^{k-2}f}{\rho}\label{eq:acceleration-flowable}\\
\frac{\partial v}{\partial t} & =-\left(\frac{\partial^{2}C}{\partial t^{2}}+v\,\frac{\partial^{2}C}{\partial s\partial t}+\frac{k+1}{t}\frac{\partial C}{\partial t}\right)\cdot\frac{\partial C}{\partial s}-\left(\frac{\partial v}{\partial s}+\frac{k+1}{t}\right)v\nonumber 
\end{align}
Notice in this flowable conserved mass model, that only the normal
component of the curve acceleration $\frac{\partial^{2}C}{\partial t^{2}}$
is governed by the Euler-Lagrange equation. The tangential acceleration,
even though it affects the internal mass flow, can be chosen freely.
We may exploit this degree of freedom to keep the tangential velocity
of the curve equal to zero, thus keeping the evolution purely geometric.
Accordingly, and just as in the constant density case, we may convert
the second order system (\ref{eq:acceleration-flowable}) into a first
order geometric system of PDE's. In particular, if we start out with
zero initial velocity, we obtain the following equivalent system of
three coupled evolution PDE's for $C$, $V$, and $\rho$ which, in
contrast with the constant density scheme, also avoids the calculation
of curvature.
\begin{align}
\overbrace{\frac{\partial V}{\partial t}}^{\mbox{acceleration}}\!\!\!\!\!\!\!\! & =\frac{\lambda k^{2}t^{k-2}}{\rho}\overbrace{f\,N}^{-\mbox{gradient}}-\overbrace{\left(V\cdot\frac{\partial C}{\partial s}\right)\frac{\partial V}{\partial s}}^{\mbox{advection}}-\overbrace{\frac{k+1}{t}V}^{\mbox{friction}}\label{eq:accelerated-system-flowable}\\
\frac{\partial C}{\partial t}= & \underbrace{(V\cdot N)}_{\beta}N,\qquad\underbrace{\frac{\partial\rho}{\partial t}=-\left(V\cdot\frac{\partial C}{\partial s}\right)\frac{\partial\rho}{\partial s}-\rho\,\frac{\partial V}{\partial s}\cdot\frac{\partial C}{\partial s}}_{\mbox{mass preservation}}\nonumber 
\end{align}
Here the velocity field defined as $V=vT+\beta N$ captures both the
tangential flow of the mass as well as the normal flow of the curve
itself. As in the constant density model, we see that the evolution
of the contour remains purely geometric (only in the normal direction),
and thus with suitable spatial extension functions $\hat{V}(x,t)$
and $\hat{\rho}(x,t)$ this system can easily be adapted to the level
set framework as well. One notable difference, however, is that evolution
equation for the level set function $\psi$ itself, becomes linear
in this case.
\begin{equation}
\frac{\partial\psi}{\partial t}=\hat{V}\cdot\nabla\psi\label{eq:levelset-flowable}
\end{equation}

\subsection{Mixing optimization dynamics with physical time dynamics\label{sec:physical}}

Again, while the kinetic energy models, including their attributed
mass density functions $\rho$, are invented purely for the sake of
improved optimization, there may be applications in which physically
meaningful considerations could nonetheless be usefully blended into
the optimization dynamics. Two particular application areas where
very strong connections could be made include dynamic tracking as
well as optimal mass transport.

\subsubsection{Connections with dynamic tracking}

Niethammer \emph{et. al.} \cite{Niethammer03,Niethammer04,Niethammer06}
introduced a new geometric dynamical active contour model that has
strong connections to the present work. The motivation is visual tracking.
The authors make the point that the use of active contours is typically
preformed statically. More specifically, the active contour captures
the given object at a certain time $t$ and then some prediction procedure
is employed to give a reasonable initial placement at time $t+1$.
The problem is that the curve evolution gets decoupled from the dynamics
of the target. The standard dynamic approaches are marker particle
based and thus lose the advantages of the level set methodology. the
shortcomings of such particle-based implementations. The works of
\cite{Niethammer03,Niethammer04,Niethammer06} develop a straightforward,
efficient, level set based approach for dynamic curve evolution which
removes the separation of segmentation and prediction, while preserving
the many advantages of level set formulations. The key idea is based
on the minimization of novel energy functional that adds dynamics
into the geodesic active contour framework.

More precisely, the above approach develops dynamical geodesic snake
models for visual tracking based on the classical action ${\bf T}-{\bf U}$
using constant density mass models. This endowed the moving contour
with dynamics in actual \emph{physical time} which could be used in
the context of dynamic observers \cite{Niethammer08}. 

Such a scheme for frame-to-fame evolution of a contour within a video
would pair very naturally, for example, with the simplest-case optimization
dynamics from Section \ref{sec:constant} using the same kinetic energy
model (\ref{eq:kinetic-constant}), but in the context of the generalized
action (\ref{eq:time-action}) for static optimization within each
individual video frame. Conversely, the more general kinetic energy
models outlined in Section \ref{sec:conserved} for optimization using
the generalized action, could be similarly be adapted to the the problem
of visual tracking using the classical action.

\subsubsection{Connections with optimal mass transport}

The conserved evolutionary mass model underlying the accelerated system
(\ref{eq:accelerated-system-flowable}) begins to exhibit clear connections
to problems in optimal mass transport \cite{Kantorovich48,Benamou00,Angenent03,Villani03},
especially in the fluid-dynamical formulation of Benamou and Brenier
\cite{Benamou00}. 

Optimal mass transport is a very old problem first introduced by the
civil engineer Monge in 1781 \cite{Monge} and concerned finding the
optimal way, in the sense of minimal transportation cost, of moving
a pile of soil from one site to another. This problem of \emph{optimal
mass transport} (OMT) was given a modern formulation in the work of
Kantorovich \cite{Kantorovich42,Kantorovich48}, and so is now known
as the Monge-{}-Kantorovich (MK) problem. As originally formulated,
the problem has no explicit dynamics, and basically leads to a metric
on probability densities, the \emph{Wasserstein distance}. Optimal
mass transport is a very active area of research with applications
to numerous disciplines including probability, econometrics, fluid
dynamics, automatic control, transportation, statistical physics,
shape optimization, expert systems, and meteorology. 

A major development in optimal mass transport theory was realized
in the seminal dynamic approach to optimal mass transport by Benamou
and Brenier \cite{Benamou00}. These authors base their approach to
OMT on ideas from fluid mechanics via the minimization of a kinetic
energy functional subject to a continuity constraint. 

The work described above is very much in line with the latter dynamics
approach. In fact, given that the mass is introduced as an independent
auxiliary variable for the sake of acceleration, we may just as easily
allow it to live within the contour interior rather than along the
contour boundary. The resulting region based extension of the kinetic
energy model (\ref{eq:kinetic-flowable}) would then match the functional
whose minimizer, as demonstrated Benamou and Brenier, yields a flow
of diffeormorphisms which minimize the Wasserstein distance between
the mass distributions at any two instances along its trajectory (including
the initial and final distributions).

\subsection{Accelerated Active Surfaces\label{sec:surfaces}}

The accelerated active contour models developed in Section \ref{sec:contours}
offer a more robust evolution framework for generic contour based
optimization problems, just as the class of Sobolev active contour
models introduced earlier. Both methodologies regularize the optimization
process, without imposing regularity on the final optimized result,
greatly boosting the evolving contour's resistance to spurious or
shallow local minimizers. In both cases, this desirable property is
achieved by effectively averaging contributions from several local
gradient forces in order to determine the instantaneous evolution
of any given point on the curve. 

In the case of Sobolev active contours, this averaging is done spatially
at each fixed time instant by an effective convolution along the curve.
Unfortunately, while special tricks exist to do this quickly for closed
curves, they do not apply to surfaces or higher dimensional manifolds,
where Laplace-Beltrami style PDE's must instead be solved along the
surface at every time instant in order to calculate the Sobolev gradeint. 

Accelerated active contour models, on the other hand, perform a temporal
rather than spatial averaging. As a particle along the curve accelerates,
its instantaneous velocity represents the accumulation of local gradient
information over its recently traveled trajectory, rather than the
accumulation of local gradient information from its neighboring contour
points at the same instant in time. An important advantage of the
time-based averaging, in contrast with the instantaneous spatially-based
averaging\footnote{In Section \ref{sec:spatial} we show how additional strategies within
the accelerated framework can be devised to further incorporate some
level of spatial averaging, thereby obtain the maximum amount of evolution
robustness and leveraging the best of both Sobolev and accelerated
optimization yet without the added computational cost of inverting
the Sobolev operator.} in Sobolev style active contours, is that the same computational
speed up in 2D will apply equally in 3D and higher. 

In the case of geometric active surfaces, we start with a potential
energy which depends only upon the geometric shape of the contour
$S$ (again, as in the contour case, not its parameterization). Under
these assumptions the first variation of the potential energy will
have the following form
\begin{align*}
\delta{\bf U} & =-\int_{S}f\,(\delta C\cdot N)\,dA
\end{align*}
where $fN$ represents a force along the unit normal $N$ at each
point on the surface $S$ and where $dA$ denotes the surface area
measure. The implicit level set framework is particularly convenient
for active surfaces given the complexities of dealing with 3D meshes.
In the level set framework, the (non-accelerated) gradient descent
surface evolution PDE has the same form as in 2D, but is applied to
a 3D grid instead. Namely
\[
\frac{\partial\psi}{\partial t}=-\hat{f}\|\nabla\psi\|
\]
where $\hat{f}(x,t)$ denotes a spatial extension of $f$ to points
away from the surface. Narrow band methods are especially important
in 3D to keep the computational cost of updating the level set function
$\psi$ to a minimum (as well as limiting the neighborhood where extension
functions such as $\hat{f}$ need to be computed and evolved).

In the simplest constant density model case, applied to surfaces.
the kinetic energy term for the accelerated model will have a similar
form but with the density $\rho$ interpreted per unit surface area.

\[
{\bf T}=\frac{1}{2}\rho\int_{S}\left(\frac{\partial S}{\partial t}\cdot\frac{\partial S}{\partial t}\right)dA
\]
 Computing the Euler-Lagrange equation of the generalized action integral
(\ref{eq:action-integral}) and writing it in the level set framework
yields the same system of first order PDE's as in the contour case,
except now in 3D dimensions,
\[
\frac{\partial\psi}{\partial t}=-\hat{\beta}\|\nabla\psi\|,\qquad\frac{\partial\hat{\beta}}{\partial t}=\frac{\lambda k^{2}t^{(k-2)}}{\rho}\hat{f}+\nabla\cdot\left(\frac{1}{2}\hat{\beta}^{2}\frac{\nabla\psi}{\|\nabla\psi\|}\right)-\frac{k+1}{t}\hat{\beta}
\]
 where $\hat{f}(x,t)$ and $\hat{\beta}(x,t)$ denote 3D spatial extensions
of $f$ and $\beta$ respectively.

\subsection{Acceleration with spatial regularity (capturing Sobolev gradient
properties)\label{sec:spatial}}

There are several ways in the PDE framework that we may seek to combine
the spatial averaging of gradient information inherent to Sobolev
gradient descent with the temporal averaging of gradient information
inherent to acceleration, while still remaining fully within the accelerated
framework, bypassing the linear operator inversion required in the
Sobolev framework. We present two different strategies for obtaining
the best-of-both. 

\subsubsection{Adding velocity diffusion}

A simple way to incorporate spatial averaging in the acceleration
process would be to heuristically add a diffusion term in the velocity
update. For a concrete example, in the conserved flowable-mass acceleration
strategy for active contours outlined in Section \ref{sec:conserved},
we could augment the acceleration PDE (\ref{eq:accelerated-system-flowable})
as follows (the coupled density evolution PDE would remain the same)
\begin{align*}
\overbrace{\frac{\partial V}{\partial t}}^{\mbox{acceleration}} & =\frac{\lambda k^{2}t^{k-2}}{\rho}\overbrace{f\,N}^{\mbox{gradient}}-\overbrace{\left(V\cdot\frac{\partial C}{\partial s}\right)\frac{\partial V}{\partial s}}^{\mbox{advection}}-\overbrace{\frac{k+1}{t}V}^{\mbox{friction}}+\overbrace{\tau\frac{\partial^{2}V}{\partial s^{2}}}^{\mbox{diffusion}}
\end{align*}
where $\tau>0$ represents a tunable diffusion coefficient. Large
values of $\tau$ would give preferential treatment to coarse scale
deformations of the evolving contour during the early stages of evolution,
with finer scale deformations gradually folding in more and more as
the contour converges toward a steady state configuration.

Such a coarse-to-fine behavior would be consistent with that of a
Sobolev active contour. In fact, diffusion over a finite amount of
time is similar to convolution with a smoothing kernel, which is indeed
one way to relate the velocity field of a Sobolev active contour with
the simple gradient field $fN$. As such, the incorporation of a diffusion
term into the acceleration PDE is the closest and most direct way
to endow the accelerated active contour with additional coarse-to-fine
Sobolev active contour behaviors without directly employing Sobolev
norms in the definition of the kinetic energy (which would would require
full linear operator inversion at every time step during the accelerated
flow, just as in actual Sobolev gradient flows).

A key difference of such an added diffusion term, compared to Sobolev
active contours, is that this smoothing process of the gradient field
along the contour is carried out concurrently with the accelerated
contour evolution itself, rather than statically at each separate
time step. As such, if the diffusion coefficient $\tau$ is small
enough to allow stable discretization of the PDE with the same time
step dictated by the other first order terms, then no additional computational
cost is incurred. As the diffusion coefficient is increased, however,
the discrete CFL conditions arising from the added second-order diffusion
term will begin to dominate in the numerical implementation of the
PDE and require smaller and smaller time steps. This could significantly
increase the computational cost as more and more numerical iterations
will be needed to simulate the same amount of accelerated flow time. 

Given that a sufficiently small amount of diffusion costs essentially
nothing in the PDE discretization, however, it doesn't make sense
to ignore this benefit from an optimization standpoint. Methodical
schemes guided purely by numerical considerations can be devised to
add velocity diffusion coefficients that will maximally boost the
regularity of the accelerated evolution with minimal or no added computational
cost. Such \emph{free gains} from small amounts of diffusion may be
stretched the farthest by allowing variable diffusion coefficients
which can be chosen based on evolving CFL conditions relevant to the
PDE discreteizations prior to considering the added diffusion terms.

\subsubsection{Incorporating mass potential energy\label{sec:mass-potential}}

An independent approach that would add spatial regularization to the
acceleration process, again without imposing any added regularity
to the converged result, would be to attach not only a kinetic energy
term to the artificially attributed mass, but also an extra potential
energy term ${\bf U}_{mass}$ which favors a smoother evolution of
the mass itself (and therefore of the object to which the auxiliary
density function is attributed). This opens up a whole new design
feature for accelerated PDE's which would allow us to incorporate
coarse-to-fine evolution properties which are qualitatively similar
to those of Sobolev gradient flows, but without the heavy computational
cost.

We foresee at least two criteria that should be satisfied when designing
the mass potential energy term ${\bf U}_{mass}$
\begin{enumerate}
\item The minimum achievable mass potential energy should be independent
of the configuration of the original variable being optimized (for
example, in the active contour case, it should achievable for any
possible contour shape) so that the final converged result, which
will correspond to a locally minimal total potential energy, will
not be influenced by the added mass potential energy term but only
by the original potential energy term to be minimized. As such, the
incorporation of ${\bf U}_{mass}$ will affect only the accelerated
evolution dynamics, without changing the original energy landscape.
\item The first variation $\delta{\bf U}_{mass}$ should not contain second
or higher order derivatives of the density function $\rho$ (nor of
its flow velocity $V$) which would, like the diffusion strategy described
earlier, impose stronger discrete time step restrictions on the numerical
discretization of the accelerated PDE system.
\end{enumerate}
In order to work out a concrete example, we revisit the accelerated
active contour model using the conserved flowable mass strategy outlined
in Section \ref{sec:conserved}, in which we suggested that the evolving
contour may be thought of as a moving container of fluid (the attributed
mass variable), and that the fluid is pushed around by the moving
container while also flowing within the container. If the fluid is
compressible, then its density can vary during this evolution, otherwise
it must remain constant, which undermines the flexibility of this
scheme compared to the simpler constant density scheme already developed
beforehand in Section \ref{sec:constant}. Yet we can still give physical
intuition to the more flexible  flowable mass model, even if we consider
the mass as an incompressible fluid. We simply imagine that the fluid
has a variable height at each point within its container (in this
case, along the contour). This allows us to naturally define a potential
energy for the mass configuration, by relating the density function
$\rho$ to the fluid height.

Using this fluid height model, we may construct the mass potential
energy connected with an arclength increment $ds$ along the curve
by first noting that the associated mass differential is given by
$dm=\rho\,ds$ and then equating the mass density $\rho$ along the
contour to a constant fluid density $\sigma$ scaled by the local
fluid height $h$. Given that the average height of the fluid column
over $ds$ would be $h/2$, we may write its potential energy as $\frac{h}{2}\,g\,dm$
where $g$ represents a gravitational constant. Combining these relationships
yields
\[
d{\bf U}=\frac{g}{2\sigma}\,\rho^{2}ds
\]
which, if we choose $\sigma=1$ (without any loss of generality since
$g$ can be chosen arbitrarily), gives the following expression for
the mass potential energy.
\[
g\int_{C}\frac{1}{2}\rho^{2}ds
\]
However, while this satisfies our second criterion (its first variation
will not involve second order derivatives of $\rho$ or of the flow
velocity $V$), it fails our first criterion. To see this, note that
lowest potential energy mass distribution for a given curve (subject
to the conservation constraint) is achieved by the constant height
distribution $\rho=\frac{M}{L}$, where $M$ denotes the constant
total conserved mass and where $L=\int_{C}ds$ denotes the total arclength
of the contour. 
\[
\min_{\rho}\left(g\int_{C}\frac{1}{2}\rho^{2}ds\right)=\frac{g}{2}\frac{M^{2}}{L}
\]
From this expression, we can see that scaling this potential energy
by the length of the curve will make the minimum achievable potential
energy $\frac{g}{2}M^{2}$ become independent of the curve $C$. This
leads to the following candidate for a mass potential energy which
also satisfies our first criterion.
\begin{equation}
{\bf U}_{mass}=g\,L\int_{C}\frac{1}{2}\rho^{2}ds\label{eq:mass-potential}
\end{equation}

Adding this to the purely contour based potential energy (Section
\ref{sec:potential}), which does not depend on the artificially added
mass, and recomputing the Euler-Lagrange equations for the generalized
action integral (\ref{eq:action-integral}) will yield a new system
of accelerated PDE's in which the gradient forces influencing the
acceleration will depend both on the mass distribution as well as
the functional to be minimized. Since the minimum constant density
mass potential can be achieved for any contour configuration, we know
that at steady state, we will have a constant mass density. If we
initialize with a constant mass density as well, then the acceleration
dynamics will favor (but not constrain) moving the mass along evolution
paths that keep the density spatially constant. Translations or uniform
rescaling of the curve would therefore become preferential evolutions,
just as for Sobolev active contours, especially with larger choices
of the tunable gravitational constant $g$.

\subsection{Incorporating stochastic acceleration terms}

Finally, the accelerated PDE framework, unlike the gradient descent
PDE framework, offers a numerical opportunity to introduce random
noise into the evolution process without destroying the continuity
of the evolution process nor of the evolving object. For example,
in the active contour acceleration scheme (\ref{eq:accelerated-system-flowable}),
we could replace the added diffusion term suggested in Section \ref{sec:spatial},
with a stochastic term as follows

\[
\overbrace{\frac{\partial V}{\partial t}}^{\mbox{acceleration}}\!\!\!\!\!\!\!\!=\frac{\lambda k^{2}t^{k-2}}{\rho}\overbrace{f\,N}^{\mbox{gradient}}-\overbrace{\left(V\cdot\frac{\partial C}{\partial s}\right)\frac{\partial V}{\partial s}}^{\mbox{advection}}-\overbrace{\frac{k+1}{t}V}^{\mbox{friction}}+\overbrace{\tau{\cal W}}^{\mbox{noise}},\qquad\frac{\partial C}{\partial t}=(V\cdot N)N
\]
where ${\cal W}$ represents samples drawn from a random noise process
and $\tau$ is a positive tunable coefficient (similar to the diffusion
coefficient in Section \ref{sec:spatial}). Since the noise is added
to the acceleration, it gets twice integrated in the construction
of the updated contour (or surface) and therefore does not immediately
interfere with the continuity nor the first order differentiability
of the evolving variable. As such, both the velocity $V$ itself as
well as the unit normal $N$ of the contour, remain continuous for
the the coupled contour evolution equation. The contour therefore
maintains regularity (at least short term). Furthermore, since upwind
differencing methods are utilized in the numerical calculation of
$\frac{\partial V}{\partial s}$ in the acceleration advection term,
discontinuities in the first derivative of $V$ do not pose a problem
as only one-sided derivatives are required. In the case of shocks,
a viscosity solution will be approximated by a proper discretization.

Adding random noise to a standard (non-accelerated) gradient descent
contour PDE, on the other hand,
\[
\overbrace{\frac{\partial C}{\partial t}}^{\mbox{velocity}}=\overbrace{f\,N}^{\mbox{gradient}}+\overbrace{\tau{\cal W}}^{\mbox{noise}}
\]
has never been a viable option since noise added directly to the velocity
is integrated only once, which does not maintain continuity in the
unit normal $N$ of the evolving contour. As such, the contour would
immediately become irregular. As such, accelerated PDE's open up a
whole new avenue for the inclusion of stochastic terms (as often exploited
in finite dimensional optimization problems) which offer an additional
strategy for increased resistance to spurious or shallow local minimizers.
The potential benefit of such a random noise term would be to provide
a second and independent mechanism (beyond the acceleration itself)
to perturb the optimization flow away from saddle points or shallow
minimizers. Once kinetic energy has been accumulated, the added benefit
of such a term is likely to be negligible. However, unlucky initializations
(assuming zero initial velocity) near local minimizers or saddle points,
could benefit from a noise driven term in the early stages while momentum
is just begining to accumulate. Note that such a strategy is not the
same as stochastic gradient descent, and should not be confused with
the recently developed PDE methods in \cite{Chaudhari17deep,Chaudhari17stochastic}
which specifically improve upon stochastic gradinent descent technqiues
used in training deep neural networks.

\section{Appendix: Derivations of various numbered equations}

\noindent\fbox{\begin{minipage}[t]{1\columnwidth - 2\fboxsep - 2\fboxrule}%
\begin{center}
\noun{\scriptsize{}calculation of equation (\ref{eq:frame-evolution})}
\par\end{center}{\scriptsize \par}
\vspace{-3mm}

{\scriptsize{}Differentiating (\ref{eq:velocity}) with respect to
the arclength parameter $s$ yields
\begin{align*}
\frac{\partial^{2}C}{\partial s\partial t} & =\frac{\partial\alpha}{\partial s}T+\alpha\underbrace{\frac{\partial T}{\partial s}}_{\kappa N}+\frac{\partial\beta}{\partial s}N+\beta\underbrace{\frac{\partial N}{\partial s}}_{-\kappa T}=\left(\frac{\partial\alpha}{\partial s}-\beta\kappa\right)T+\left(\frac{\partial\beta}{\partial s}+\alpha\kappa\right)N
\end{align*}
and differentiating $T=\frac{\partial C}{\partial s}$ yields}{\scriptsize \par}

{\scriptsize{}
\begin{align*}
\frac{\partial T}{\partial t} & =\frac{\partial}{\partial t}\frac{\partial C}{\partial s}=\frac{\partial}{\partial t}\left(\frac{\frac{\partial C}{\partial p}}{\left\Vert \frac{\partial C}{\partial p}\right\Vert }\right)=\frac{\frac{\partial^{2}C}{\partial t\partial p}}{\left\Vert \frac{\partial C}{\partial p}\right\Vert }-\frac{\frac{\partial C}{\partial p}}{\left\Vert \frac{\partial C}{\partial p}\right\Vert ^{2}}\frac{\partial}{\partial t}\left\Vert \frac{\partial C}{\partial p}\right\Vert =\frac{\frac{\partial^{2}C}{\partial p\partial t}}{\left\Vert \frac{\partial C}{\partial p}\right\Vert }-\frac{\frac{\partial C}{\partial p}}{\left\Vert \frac{\partial C}{\partial p}\right\Vert ^{2}}\frac{\frac{\partial^{2}C}{\partial p\partial t}\cdot\frac{\partial C}{\partial p}}{\left\Vert \frac{\partial C}{\partial p}\right\Vert }\\
 & =\frac{\partial^{2}C}{\partial s\partial t}-T\left(\frac{\partial^{2}C}{\partial s\partial t}\cdot T\right)=\left(\frac{\partial^{2}C}{\partial s\partial t}\cdot N\right)N=\left(\frac{\partial\beta}{\partial s}+\alpha\kappa\right)N
\end{align*}
which gives the first part of (\ref{eq:frame-evolution}) with the
second part due the rotation relationship between $T$ and $N$.}{\scriptsize \par}%
\end{minipage}}\bigskip{}
\noindent\fbox{\begin{minipage}[t]{1\columnwidth - 2\fboxsep - 2\fboxrule}%
\begin{center}
\noun{\scriptsize{}calculation of equation (\ref{eq:acceleration-constant})}
\par\end{center}{\scriptsize \par}
\vspace{-3mm}

{\scriptsize{}Letting $C(p,t)$ denote a parameterization of the evolving
curve $C$ with a time-independent spatial parameter $p$ and with
$s$ denoting the time-dependent arclength parameter we compute (ignoring
temporary boundary terms when applying integration by parts and assuming
a closed curve so that spatial boundary terms cancel):
\begin{align*}
 & \delta\int_{0}^{1}\frac{t^{k+1}}{k}\left({\bf T}-\lambda k^{2}t^{k-2}{\bf U}\right)\,dt=\delta\int_{0}^{1}\left(\int_{0}^{1}\frac{1}{2}\frac{t^{k+1}}{k}\rho\frac{\partial C}{\partial t}\cdot\frac{\partial C}{\partial t}\,ds-\lambda kt^{2k-1}{\bf U}\right)\,dt=\int_{0}^{1}\delta\left(\int_{0}^{1}\frac{1}{2}\frac{t^{k+1}}{k}\rho\frac{\partial C}{\partial t}\cdot\frac{\partial C}{\partial t}\left\Vert \frac{\partial C}{\partial p}\right\Vert dp-\lambda kt^{2k-1}{\bf U}\right)\,dt\\
 & =\int_{0}^{1}\frac{\rho}{k}\left(\int_{0}^{1}\frac{1}{2}t^{k+1}\left(\frac{\partial C}{\partial t}\cdot\delta\frac{\partial C}{\partial t}\left\Vert \frac{\partial C}{\partial p}\right\Vert +\frac{1}{2}\left\Vert \frac{\partial C}{\partial t}\right\Vert ^{2}\delta\left\Vert \frac{\partial C}{\partial p}\right\Vert \right)dp-\frac{\lambda k^{2}t^{2k-1}}{\rho}\delta{\bf U}\right)\,dt\\
 & =\int_{0}^{1}\frac{\rho}{k}\left(\int_{0}^{1}-\frac{\partial}{\partial t}\left(t^{k+1}\frac{\partial C}{\partial t}\left\Vert \frac{\partial C}{\partial p}\right\Vert \right)\cdot\delta C+\frac{1}{2}t^{k+1}\left\Vert \frac{\partial C}{\partial t}\right\Vert ^{2}\delta\frac{\partial C}{\partial p}\cdot\frac{\partial C}{\partial s}\,dp+\frac{\lambda k^{2}t^{2k-1}}{\rho}\int_{C}f\,(\delta C\cdot N)\,ds\right)\,dt\\
 & =\int_{0}^{1}\frac{\rho}{k}\left(\int_{0}^{1}-t^{k+1}\frac{\partial^{2}C}{\partial t^{2}}\left\Vert \frac{\partial C}{\partial p}\right\Vert -(k+1)t^{k}\frac{\partial C}{\partial t}\left\Vert \frac{\partial C}{\partial p}\right\Vert -t^{k+1}\frac{\partial C}{\partial t}\frac{\partial}{\partial t}\left\Vert \frac{\partial C}{\partial p}\right\Vert -\frac{1}{2}t^{k+1}\frac{\partial}{\partial p}\left(\left\Vert \frac{\partial C}{\partial t}\right\Vert ^{2}\frac{\partial C}{\partial s}\right)\cdot\delta C\,dp+\frac{\lambda k^{2}t^{2k-1}}{\rho}\int_{C}fN\cdot\delta C\,ds\right)\,dt\\
 & =\int_{0}^{1}\frac{t^{k+1}}{k}\rho\left(\int_{0}^{1}\left(-\frac{\partial^{2}C}{\partial t^{2}}\left\Vert \frac{\partial C}{\partial p}\right\Vert -\frac{k+1}{t}\frac{\partial C}{\partial t}\left\Vert \frac{\partial C}{\partial p}\right\Vert -\frac{\partial C}{\partial t}\left(\frac{\partial^{2}C}{\partial p\partial t}\cdot\frac{\partial C}{\partial s}\right)-\frac{1}{2}\frac{\partial}{\partial p}\left(\left\Vert \frac{\partial C}{\partial t}\right\Vert ^{2}\frac{\partial C}{\partial s}\right)\right)\cdot\delta C\,dp+\frac{\lambda k^{2}t^{k-2}}{\rho}\int_{C}fN\cdot\delta C\,ds\right)\,dt\\
 & =\int_{0}^{1}\frac{t^{k+1}}{k}\rho\int_{C}\underbrace{\left(-\frac{\partial^{2}C}{\partial t^{2}}-\frac{k+1}{t}\frac{\partial C}{\partial t}-\left(\frac{\partial^{2}C}{\partial s\partial t}\cdot\frac{\partial C}{\partial s}\right)\frac{\partial C}{\partial t}-\frac{\partial}{\partial s}\left(\frac{1}{2}\left\Vert \frac{\partial C}{\partial t}\right\Vert ^{2}\frac{\partial C}{\partial s}\right)+\frac{\lambda k^{2}t^{k-2}}{\rho}fN\right)}_{\mbox{Set to zero for Euler-Lagrange equation}}\cdot\delta C\,ds\,dt
\end{align*}
}{\scriptsize \par}%
\end{minipage}}\bigskip{}
\noindent\fbox{\begin{minipage}[t]{1\columnwidth - 2\fboxsep - 2\fboxrule}%
\begin{center}
\noun{\scriptsize{}calculation of equation (\ref{eq:accelerated-system-constant})}
\par\end{center}{\scriptsize \par}
\vspace{-3mm}

{\scriptsize{}Decomposing the acceleration $C_{tt}$ into tangential
and normal components yields
\begin{align*}
\frac{\partial^{2}C}{\partial t^{2}} & =-\frac{k+1}{t}\frac{\partial C}{\partial t}-\left(\frac{\partial^{2}C}{\partial s\partial t}\cdot\frac{\partial C}{\partial s}\right)\frac{\partial C}{\partial t}-\left(\frac{\partial^{2}C}{\partial s\partial t}\cdot\frac{\partial C}{\partial t}\right)\frac{\partial C}{\partial s}-\frac{1}{2}\left\Vert \frac{\partial C}{\partial t}\right\Vert ^{2}\frac{\partial^{2}C}{\partial s^{2}}+\frac{\lambda k^{2}t^{k-2}f}{\rho}N\\
 & =-\frac{k+1}{t}\frac{\partial C}{\partial t}-\left(\frac{\partial\alpha}{\partial s}-\beta\kappa\right)\frac{\partial C}{\partial t}-\left(\alpha\frac{\partial\alpha}{\partial s}+\beta\frac{\partial\beta}{\partial s}\right)T-\frac{\alpha^{2}+\beta^{2}}{2}\kappa N+\frac{\lambda k^{2}t^{k-2}f}{\rho}N\\
\frac{\partial^{2}C}{\partial t^{2}}\cdot T & =-\left(\frac{k+1}{t}+\frac{\partial\alpha}{\partial s}-\beta\kappa\right)\alpha-\left(\alpha\frac{\partial\alpha}{\partial s}+\beta\frac{\partial\beta}{\partial s}\right)=-\left(\frac{k+1}{t}+2\frac{\partial\alpha}{\partial s}-\beta\kappa\right)\alpha-\beta\frac{\partial\beta}{\partial s}\\
\frac{\partial^{2}C}{\partial t^{2}}\cdot N & =-\left(\frac{k+1}{t}+\frac{\partial\alpha}{\partial s}-\beta\kappa\right)\beta-\frac{\alpha^{2}+\beta^{2}}{2}\kappa+\frac{\lambda k^{2}t^{k-2}f}{\rho}
\end{align*}
Now inserting these acceleration components into (\ref{eq:speed-evolution})
yields
\begin{align*}
\frac{\partial\alpha}{\partial t} & =\underbrace{-\left(\frac{k+1}{t}+2\frac{\partial\alpha}{\partial s}-\beta\kappa\right)\alpha-\beta\frac{\partial\beta}{\partial s}}_{\frac{\partial^{2}C}{\partial t^{2}}\cdot T}+\beta\left(\frac{\partial\beta}{\partial s}+\alpha\kappa\right)=\left(-\frac{k+1}{t}-2\frac{\partial\alpha}{\partial s}+2\beta\kappa\right)\alpha\\
\frac{\partial\beta}{\partial t} & =\underbrace{-\left(\frac{k+1}{t}+\frac{\partial\alpha}{\partial s}-\beta\kappa\right)\beta-\frac{\alpha^{2}+\beta^{2}}{2}\kappa+\frac{\lambda k^{2}t^{k-2}f}{\rho}}_{\frac{\partial^{2}C}{\partial t^{2}}\cdot N}-\alpha\left(\frac{\partial\beta}{\partial s}+\alpha\kappa\right)=-\frac{k+1}{t}\beta-\frac{\partial}{\partial s}(\alpha\beta)+\left(\frac{1}{2}\beta^{2}-\frac{3}{2}\alpha^{2}\right)\kappa+\frac{\lambda k^{2}t^{k-2}f}{\rho}
\end{align*}
Given zero initial velocity ($\alpha$=0 and $\beta$=0), simple inspection
shows that $\alpha$ remains zero, leading to the simplified evolution
(\ref{eq:accelerated-system-constant}).}{\scriptsize \par}%
\end{minipage}}\bigskip{}
\noindent\fbox{\begin{minipage}[t]{1\columnwidth - 2\fboxsep - 2\fboxrule}%
\begin{center}
\noun{\scriptsize{}calculation of equation (\ref{eq:levelset-constant})}
\par\end{center}{\scriptsize \par}
\vspace{-3mm}

{\scriptsize{}Assuming we represent the evolving curve $C(p,t)$ as
the zero level set of an evolving function $\psi(x,t)$ and letting
$\hat{\beta}(x,t)$ denote an evolving spatial extension of the evolving
normal speed function $\beta(p,t)$ along curve, then we have }{\scriptsize \par}

{\scriptsize{}\vspace{-2mm}
}{\scriptsize \par}

{\scriptsize{}
\[
\psi\left(C(p,t),t\right)=0\quad\mbox{and}\quad\hat{\beta}\left(C(p,t),t\right)=\beta(p,t)
\]
}{\scriptsize \par}

{\scriptsize{}\vspace{-2mm}
}{\scriptsize \par}

{\scriptsize{}Differentiating with respect to $t$ yields}{\scriptsize \par}

{\scriptsize{}\vspace{-4mm}
}{\scriptsize \par}

{\scriptsize{}
\[
\frac{\partial\psi}{\partial t}+\nabla\psi\cdot\frac{\partial C}{\partial t}=0\quad\mbox{and}\quad\frac{\partial\hat{\beta}}{\partial t}+\nabla\hat{\beta}\cdot\frac{\partial C}{\partial t}=\frac{\partial\beta}{\partial t}
\]
}{\scriptsize \par}

\vspace{-1mm}

{\scriptsize{}Extending the contour evolution $\frac{\partial C}{\partial t}=\beta N$
to other level sets as $\hat{\beta}\hat{N}$, where $\hat{N}=-\frac{\nabla\psi}{\|\nabla\psi\|}$
(noting that this convention for the extension of the inward unit
normal requires that the level set function be negative inside the
contour and positive outside), yields}{\scriptsize \par}

{\scriptsize{}\vspace{-2mm}
}{\scriptsize \par}

{\scriptsize{}
\[
\frac{\partial\psi}{\partial t}=\hat{\beta}\|\nabla\psi\|\quad\mbox{and}\quad\frac{\partial\hat{\beta}}{\partial t}=\frac{\partial\beta}{\partial t}+\nabla\hat{\beta}\cdot\frac{\hat{\beta}\nabla\psi}{\|\nabla\psi\|}
\]
}{\scriptsize \par}

\vspace{-2mm}

{\scriptsize{}which, after substitution of $\frac{\partial\beta}{\partial t}$
using (\ref{eq:accelerated-system-constant}) results in the level
set version of the system in (\ref{eq:levelset-constant}).}{\scriptsize \par}%
\end{minipage}}\bigskip{}
\noindent\fbox{\begin{minipage}[t]{1\columnwidth - 2\fboxsep - 2\fboxrule}%
\begin{center}
\noun{\scriptsize{}calculation of equation (\ref{eq:acceleration-flowable})}
\par\end{center}{\scriptsize \par}
{\scriptsize{}\vspace{-3mm}
}{\scriptsize \par}

{\scriptsize{}Let us introduce, along with the mass density $\rho$
and its internal flow speed $v$ with respect to the arclength parameter
$s$, corresponding variables for the mass density $\mu(p,t)$ and
internal flow speed $\xi(p,t)$ with respect to a time-independent
contour parameter $p$. These pairs of densities and internal flow
speeds are related to each other through the parameterization speed
$\left\Vert C_{p}\right\Vert $ of the contour as follows. }{\scriptsize \par}

{\scriptsize{}\vspace{-2mm}
}{\scriptsize \par}

{\scriptsize{}
\begin{equation}
\mu=\rho\left\Vert C_{p}\right\Vert \quad\mbox{and}\quad v=\xi\,\left\Vert C_{p}\right\Vert \quad(\mbox{with matching flux expressions }\mu\,\xi=\rho\,v)\label{eq:substitute}
\end{equation}
Differentiating with respect to $t$, yields the following relationships
between the density and flow speed evolution as well.}{\scriptsize \par}

{\scriptsize{}\vspace{-2mm}
}{\scriptsize \par}

{\scriptsize{}
\begin{equation}
\mu_{t}-\mu\,C_{ts}\cdot C_{s}=\rho_{t}\left\Vert C_{p}\right\Vert \quad\mbox{and}\quad v_{t}-v\,C_{ts}\cdot C_{s}=\xi_{t}\left\Vert C_{p}\right\Vert \label{eq:substitute-evolve}
\end{equation}
Applying these substitutions to the kinetic energy (\ref{eq:kinetic-flowable})
and continuity constraint (\ref{eq:continuity-constraint}) yields }{\scriptsize \par}

{\scriptsize{}\vspace{-2mm}
}{\scriptsize \par}

{\scriptsize{}
\[
{\bf T}=\int_{0}^{1}\frac{1}{2}\mu\left\Vert C_{t}+\xi\,C_{p}\right\Vert ^{2}dp\quad\mbox{with mass continuity constraint}\quad\mu_{t}+\left(\mu\,\xi\right)_{p}=0
\]
We plug this into the generalized action integral (\ref{eq:action-integral})
with a Lagrange multiplier function $\lambda(p,t)$ and compute the
first variation.}{\scriptsize \par}

{\scriptsize{}\vspace{-5mm}
}{\scriptsize \par}

{\scriptsize{}
\begin{align*}
\delta & \int_{0}^{1}\frac{t^{k+1}}{k}\left({\bf T}-\lambda k^{2}t^{k-2}{\bf U}\right)+\int_{0}^{1}\lambda\left(\mu_{t}+(\mu\xi)_{p}\right)dp\,dt=\int_{0}^{1}\delta\int_{0}^{1}\frac{1}{2}\frac{t^{k+1}}{k}\mu\,\left\Vert C_{t}+\xi\,C_{p}\right\Vert ^{2}+\lambda\,\left(\mu_{t}+(\mu\xi)_{p}\right)\,dp-\lambda kt^{2k-1}\delta{\bf U}\,dt\\
= & \int_{0}^{1}\int_{0}^{1}\frac{1}{2}\frac{t^{k+1}}{k}\left\Vert C_{t}+\xi\,C_{p}\right\Vert ^{2}\delta\mu+\frac{t^{k+1}}{k}\mu\,(C_{t}+\xi\,C_{p})\cdot\delta\left(C_{t}+\xi\,C_{p}\right)+\underbrace{\left(\mu_{t}+(\mu\xi)_{p}\right)}_{=0}\delta\lambda+\lambda\,\delta\left(\mu_{t}+(\mu\xi)_{p}\right)+\lambda kt^{2k-1}f\,(\delta C\cdot N)\,\|C_{p}\|\,dp\,dt\\
= & \int_{0}^{1}\int_{0}^{1}\frac{1}{2}\frac{t^{k+1}}{k}\left\Vert C_{t}+\xi\,C_{p}\right\Vert ^{2}\delta\mu+\frac{t^{k+1}}{k}\mu\,(C_{t}+\xi\,C_{p})\cdot C_{p}\,\delta\xi-\lambda_{t}\delta\mu-\lambda_{p}\delta\left(\mu\xi\right)-\left(\frac{t^{k+1}}{k}\mu\,(C_{t}+\xi\,C_{p})\right)_{t}\cdot\delta C-\left(\frac{t^{k+1}}{k}\mu\xi\,(C_{t}+\xi\,C_{p})\right)_{p}\cdot\delta C\\
 & \qquad+\lambda kt^{2k-1}f\,N\cdot\delta C\,\|C_{p}\|\,dp\,dt\\
= & \int_{0}^{1}\int_{0}^{1}\left(\frac{1}{2}\frac{t^{k+1}}{k}\left\Vert C_{t}+\xi\,C_{p}\right\Vert ^{2}-\lambda_{t}-\lambda_{p}\xi\right)\,\delta\mu+\mu\,\left(\frac{t^{k+1}}{k}(C_{t}+\xi\,C_{p})\cdot C_{p}-\lambda_{p}\right)\,\delta\xi\\
 & \qquad-\frac{t^{k+1}}{k}\mu\left(\underbrace{\left(\mu_{t}+(\mu\xi)_{p}\right)}_{=0}\,\frac{C_{t}+\xi\,C_{p}}{\mu}+\left(C_{t}+\xi\,C_{p}\right)_{t}+\xi\,\left(C_{t}+\xi\,C_{p}\right)_{p}+\frac{k+1}{t}(C_{t}+\xi\,C_{p})-\frac{\lambda k^{2}t^{k-2}f\,N}{\rho}\right)\cdot\delta C\,dp\,dt
\end{align*}
}{\scriptsize \par}

{\scriptsize{}\vspace{-2mm}
}{\scriptsize \par}

{\scriptsize{}The optimality conditions with respect to variations
$\delta\xi$ and $\delta\mu$ respectively yield}{\scriptsize \par}

{\scriptsize{}\vspace{-2mm}
}{\scriptsize \par}

{\scriptsize{}
\[
\lambda_{p}=\frac{t^{k+1}}{k}\left(\xi\,\|C_{p}\|^{2}+C_{t}\cdot C_{p}\right)\quad\mbox{and}\quad\lambda_{t}=\frac{1}{2}\frac{t^{k+1}}{k}\left\Vert C_{t}+\xi\,C_{p}\right\Vert ^{2}-\xi\lambda_{p}
\]
which, when combined, give the following evolution for the Lagrange
multiplier}{\scriptsize \par}

{\scriptsize{}\vspace{-5mm}
}{\scriptsize \par}

{\scriptsize{}
\begin{align*}
\lambda_{t} & =\frac{t^{k+1}}{k}\left(\frac{1}{2}\left\Vert C_{t}+\xi\,C_{p}\right\Vert ^{2}-\left(\|\xi\,C_{p}\|^{2}+C_{t}\cdot\xi C_{p}\right)\right)=\frac{t^{k+1}}{k}\left(\frac{1}{2}\left\Vert C_{t}+v\,C_{s}\right\Vert ^{2}-\left(\|v\,C_{s}\|^{2}+C_{t}\cdot v\,C_{s}\right)\right)=\frac{1}{2}\frac{t^{k+1}}{k}\left(\|C_{t}\|^{2}-v^{2}\right)
\end{align*}
}{\scriptsize \par}

{\scriptsize{}\vspace{-2mm}
}{\scriptsize \par}

{\scriptsize{}We eliminate the Lagrange multiplier by equating $\lambda_{tp}$
and $\lambda_{pt}$ to obtain the following internal flow speed evolution}{\scriptsize \par}

{\scriptsize{}\vspace{-5mm}
}{\scriptsize \par}

{\scriptsize{}
\begin{align*}
0= & \lambda_{tp}-\lambda_{pt}=\frac{1}{2}\frac{t^{k+1}}{k}\left(\|C_{t}\|^{2}-v^{2}\right)_{p}-\left(\frac{t^{k+1}}{k}\left(C_{t}+\xi\,C_{p}\right)\cdot C_{p}\right)_{t}\\
0= & \underbrace{\frac{t^{k+1}}{k}}_{\mbox{drop}}\left(\underbrace{C_{t}\cdot C_{tp}}_{\mbox{cancel}}-vv_{p}-C_{tt}\cdot C_{p}-\xi_{t}C_{p}\cdot C_{p}-2\xi\,C_{tp}\cdot C_{p}-\underbrace{C_{t}\cdot C_{tp}}_{\mbox{cancel}}-\frac{k+1}{t}\left(C_{t}+\xi\,C_{p}\right)\cdot C_{p}\right)\\
0= & \underbrace{\|C_{p}\|}_{\mbox{drop}}\left(-vv_{s}-C_{tt}\cdot C_{s}-\underbrace{(\xi_{t}\|C_{p}\|+v\,C_{ts}\cdot C_{s})}_{v_{t}\mbox{ by using \eqref{eq:substitute-evolve}}}-v\,C_{ts}\cdot C_{s}-\frac{k+1}{t}\left(C_{t}+v\,C_{s}\right)\cdot C_{s}\right)\\
v_{t}= & -\left(C_{tt}+v\,(C_{t}+v\,C_{s})_{s}+\frac{k+1}{t}\left(C_{t}+v\,C_{s}\right)\right)\cdot C_{s}=-\left(C_{tt}+vC_{ts}+\frac{k+1}{t}C_{t}\right)\cdot C_{s}-\left(v_{s}+\frac{k+1}{t}\right)v
\end{align*}
}{\scriptsize \par}

{\scriptsize{}\vspace{-2mm}
}{\scriptsize \par}

{\scriptsize{}Finally, the optimality condition with respect to the
curve perturbation $\delta C$ yields the following acceleration equation
for the contour}{\scriptsize \par}

{\scriptsize{}\vspace{-5mm}
}{\scriptsize \par}

{\scriptsize{}
\begin{align*}
0 & =C_{tt}+\!\!\!\underbrace{\xi_{t}\|C_{p}\|}_{v_{t}-vC_{ts}\cdot C_{s}}\!\!\!C_{s}+v\,C_{ts}+v\,\left(C_{t}+v\,C_{s}\right)_{s}+\frac{k+1}{t}(C_{t}+v\,C_{s})-\frac{\lambda k^{2}t^{k-2}f}{\rho}N\\
0 & =C_{tt}+\left(v_{t}-v\,C_{ts}\cdot C_{s}\right)\,C_{s}+v\,C_{ts}+v\,\left(C_{t}+v\,C_{s}\right)_{s}+\frac{k+1}{t}(C_{t}+v\,C_{s})-\frac{\lambda k^{2}t^{k-2}f}{\rho}N\\
0 & =\underbrace{\left(C_{tt}+v\,\left(C_{t}+v\,C_{s}\right)_{s}+\frac{k+1}{t}(C_{t}+v\,C_{s})+v\,C_{ts}\right)}_{\mbox{some vector}}-\underbrace{\left(\left(C_{tt}+v(C_{t}+v\,C_{s})_{s}+\frac{k+1}{t}\left(C_{t}+v\,C_{s}\right)+v\,C_{ts}\right)\cdot C_{s}\right)C_{s}}_{\mbox{its tangential component}}-\frac{\lambda k^{2}t^{k-2}f}{\rho}N\\
0 & =\underbrace{\left(C_{tt}+v(C_{t}+v\,C_{s})_{s}+\frac{k+1}{t}\left(C_{t}+v\,C_{s}\right)+vC_{ts}\right)\cdot N}_{\mbox{its normial projection}}-\frac{\lambda k^{2}t^{k-2}f}{\rho}=\left(C_{tt}+2v\,C_{ts}+v^{2}C_{ss}+\frac{k+1}{t}C_{t}\right)\cdot N-\frac{\lambda k^{2}t^{k-2}f}{\rho}
\end{align*}
}{\scriptsize \par}%
\end{minipage}}\bigskip{}
\noindent\fbox{\begin{minipage}[t]{1\columnwidth - 2\fboxsep - 2\fboxrule}%
\begin{center}
\noun{\scriptsize{}calculation of equation (\ref{eq:accelerated-system-flowable})}
\par\end{center}{\scriptsize \par}
{\scriptsize{}\vspace{-3mm}
}{\scriptsize \par}

{\scriptsize{}Plugging the normal and tangential (unconstrained) acceleration
components from (\ref{eq:acceleration-flowable}) into (\ref{eq:speed-evolution})
yields
\begin{align*}
\alpha_{t} & =\underbrace{C_{tt}\cdot T}_{\mbox{free}}+\beta\left(\beta_{s}+\alpha\kappa\right)\\
\beta_{t} & =\underbrace{\left(-2v\,\overbrace{\left(\beta_{s}+\alpha\kappa\right)}^{C_{ts}\cdot N}-v^{2}\overbrace{\kappa}^{C_{ss}\cdot N}-\frac{k+1}{t}\overbrace{\beta}^{C_{t}\cdot N}+\lambda k^{2}t^{k-2}\frac{f}{\rho}\right)}_{C_{tt}\cdot N\mbox{ from equation \eqref{eq:acceleration-flowable}}}-\alpha(\beta_{s}+\alpha\kappa)\\
 & =-\left(2v+\alpha\right)\beta_{s}-\left(\alpha+v\right)^{2}\kappa-\frac{k+1}{t}\beta+\lambda k^{2}t^{k-2}\frac{f}{\rho}\\
v_{t} & =-C_{tt}\cdot T-v\,\overbrace{(\alpha_{s}-\beta\kappa)}^{C_{ts}\cdot T}-\frac{k+1}{t}(\overbrace{\alpha}^{C_{t}\cdot T}+v)-vv_{s}\qquad\left[\mbox{from second part of equation \eqref{eq:acceleration-flowable}}\right]\\
\alpha_{t}+v_{t} & =\left(C_{tt}\cdot T+\beta\beta_{s}+\alpha\beta\kappa\right)+\left(-C_{tt}\cdot T-v\,(\alpha_{s}-\beta\kappa)-vv_{s}-\frac{k+1}{t}\left(\alpha+v\right)\right)\\
\left(\alpha+v\right)_{t} & =\beta\beta_{s}-\left(\alpha_{s}+v_{s}\right)v+\left(\alpha+v\right)\beta\kappa-\frac{k+1}{t}\left(\alpha+v\right)\\
\rho_{t}+(\rho v)_{s} & =-\rho\,\overbrace{\left(\alpha_{s}-\beta\kappa\right)}^{C_{ts}\cdot T}\qquad\qquad\qquad\qquad\qquad\qquad\qquad\qquad\quad\left[\mbox{from equation \eqref{eq:continuity-constraint}}\right]
\end{align*}
We can now rewrite the system as follows.
\begin{align*}
\alpha_{t} & =\underbrace{C_{tt}\cdot T}_{\mbox{free}}+\beta\left(\beta_{s}+\alpha\kappa\right)\\
\beta_{t} & =-\beta_{s}v-\left(\alpha+v\right)\beta_{s}-\left(\alpha+v\right)^{2}\kappa-\frac{k+1}{t}\beta+\lambda k^{2}t^{k-2}\frac{f}{\rho}\\
\left(\alpha+v\right)_{t} & =-\left(\alpha+v\right)_{s}v+\beta\beta_{s}+\left(\alpha+v\right)\beta\kappa-\frac{k+1}{t}\left(\alpha+v\right)\\
\rho_{t} & =-\rho_{s}v-\rho\left(\alpha+v\right)_{s}+\rho\beta\kappa
\end{align*}
where we can see that freedom to choose $C_{tt}\cdot T$ is equivalent
to freedom to choose the evolution of $\alpha$. As such, we may conveniently
choose $\alpha_{t}=0$. Assuming that we start out with zero initial
velocity ($\alpha=\beta=0$) this would mean $\alpha$remains zero,
yielding the following simplified system.
\begin{align*}
\beta_{t} & +v\beta_{s}=-v\left(\beta_{s}+v\kappa\right)-\frac{k+1}{t}\beta+\lambda k^{2}t^{k-2}\frac{f}{\rho}\\
v_{t} & +vv_{s}=\beta\left(\beta_{s}+v\kappa\right)-\frac{k+1}{t}v\\
\rho_{t} & +v\rho_{s}=\rho\left(\beta\kappa-v_{s}\right)
\end{align*}
Finally, we may transform the system by defining $V=vT+\beta N$ to
avoid the explicit calculation of curvature. Noting that
\[
V_{s}=\left(vT+\beta N\right)_{s}=\left(v_{s}-\beta\kappa\right)T+\left(\beta_{s}+v\kappa\right)N
\]
and, substituting $\alpha=0$ into (\ref{eq:frame-evolution}), to
obtain
\[
T_{t}=\beta_{s}N\quad\mbox{and}\quad N_{t}=-\beta_{s}T
\]
we may compute
\begin{align*}
V_{t} & =\left(vT+\beta N\right)_{t}=\left(v_{t}-\beta\beta_{s}\right)T+\left(\beta_{t}+v\beta_{s}\right)N\\
 & =\left(\underbrace{\beta(\beta_{s}+v\kappa)-vv_{s}-\frac{k+1}{t}v}_{v_{t}}-\beta\beta_{s}\right)T+\left(\underbrace{-v(\beta_{s}+v\kappa)-v\beta_{s}-\frac{k+1}{t}\beta+\lambda k^{2}t^{k-2}\frac{f}{\rho}}_{\beta_{t}}+v\beta_{s}\right)N\\
 & =-v\underbrace{\left(\left(v_{s}-\beta\kappa\right)T+\left(\beta_{s}+v\kappa\right)N\right)}_{V_{s}}-\frac{k+1}{t}\underbrace{\left(vT+\beta N\right)}_{V}+\left(\lambda k^{2}t^{k-2}\frac{f}{\rho}\right)N
\end{align*}
as well as
\[
\rho_{t}+\underbrace{v}_{V\cdot T}\rho_{s}=\rho\,\underbrace{(\beta\kappa-v_{s})}_{V_{s}\cdot T}
\]
}{\scriptsize \par}%
\end{minipage}}\bigskip{}

\bibliographystyle{plain}
\bibliography{refs}

\begin{thebibliography}{10}

\bibitem{Angenent03}
Sigurd Angenent, Steven Haker, and Allen Tannenbaum.
\newblock Minimizing flows for the monge-kantorovich problem.
\newblock {\em SIAM J. Math. Analysis}, 35:61--97, 2003.

\bibitem{Bardelli}
E.~Bardelli, M.~Colombo, A.~Mennucci, and A.~Yezzi.
\newblock Multiple object tracking via prediction and filtering with a
  sobolev-type metric on curves.
\newblock In {\em European Conf. Computer Vision}, pages 143--152, 2012.

\bibitem{Benamou00}
J.D. Benamou and Y.~Brenier.
\newblock A computational fluid mechanics solution to the monge-kantorovich
  mass transfer problem.
\newblock {\em Numerische Mathematik}, 84:375--393, 2000.

\bibitem{Boyd04}
Stephen Boyd and Lieven Vandenberghe.
\newblock {\em Convex optimization}.
\newblock Cambridge university press, 2004.

\bibitem{Bubeck15}
S{\'e}bastien Bubeck, Yin~Tat Lee, and Mohit Singh.
\newblock A geometric alternative to nesterov's accelerated gradient descent.
\newblock {\em CoRR}, abs/1506.08187, 2015.

\bibitem{Caselles98}
Vincent Caselles, Ron Kimmel, and Guillermo Sapiro.
\newblock Geodesic active contours.
\newblock {\em International Journal on Comptuer Vision}, 22(1):61--79, 1997.

\bibitem{Chambolle11}
Antonin Chambolle and Thomas Pock.
\newblock A first-order primal-dual algorithm for convex problems with
  applications to imaging.
\newblock {\em Journal of mathematical imaging and vision}, 40(1):120--145,
  2011.

\bibitem{Chan01}
Tony Chan and Luminita Vese.
\newblock Active contours without edges.
\newblock {\em IEEE Transactions on Image Processing}, 10(2):266--277, 2001.

\bibitem{Chan06}
Tony~F Chan, Selim Esedoglu, and Mila Nikolova.
\newblock Algorithms for finding global minimizers of image segmentation and
  denoising models.
\newblock {\em SIAM journal on applied mathematics}, 66(5):1632--1648, 2006.

\bibitem{Charpiat2005}
G.~Charpiat, R.~Keriven, J.P. Pons, and O.~Faugeras.
\newblock Designing spatially coherent minimizing flows for variational
  problems based on active contours.
\newblock In {\em Int. Conference Computer Vision}, 2005.

\bibitem{Charpiat2005_Shape_metrics}
Guillaume Charpiat, Olivier Faugeras, and Renaud Keriven.
\newblock Approximations of shape metrics and application to shape warping and
  empirical shape statistics.
\newblock {\em Foundations of Computational Mathematics}, 5(1):1--58, 2005.

\bibitem{Chaudhari17deep}
Pratik Chaudhari, Adam Oberman, Stanley Osher, Stefano Soatto, and Guillame
  Carlier.
\newblock Deep relaxation: partial differential equations for optimizing deep
  neural networks.
\newblock {\em arXiv preprint arXiv:1704.04932}, 2017.

\bibitem{Chaudhari17stochastic}
Pratik Chaudhari and Stefano Soatto.
\newblock Stochastic gradient descent performs variational inference, converges
  to limit cycles for deep networks.
\newblock {\em arXiv preprint arXiv:1710.11029}, 2017.

\bibitem{Flammarion15}
Nicolas Flammarion and Francis Bach.
\newblock From averaging to acceleration, there is only a step-size.
\newblock In {\em Proceedings of Machine Learning Research}, volume~40, pages
  658--695, 2015.

\bibitem{Ghadimi16}
Saeed Ghadimi and Guanghui Lan.
\newblock Accelerated gradient methods for nonconvex nonlinear and stochastic
  programming.
\newblock {\em Math. Program.}, 156(1-2):59--99, 2016.

\bibitem{Goldstein02}
Herbert Goldstein, Charles~P. Poole, and John~L. Safko.
\newblock {\em Classical {M}echanics}.
\newblock Addison Wesley, 2002.

\bibitem{Goldstein10}
Tom Goldstein, Xavier Bresson, and Stanley Osher.
\newblock Geometric applications of the split bregman method: segmentation and
  surface reconstruction.
\newblock {\em Journal of Scientific Computing}, 45(1-3):272--293, 2010.

\bibitem{Horn81}
Berthold~K.P. Horn and Brian~G. Schunck.
\newblock Determining optical flow.
\newblock {\em Artifical Intelligence}, 17:185--203, 1981.

\bibitem{Hu09}
Chonghai Hu, Weike Pan, and James~T. Kwok.
\newblock Accelerated gradient methods for stochastic optimization and online
  learning.
\newblock In Y.~Bengio, D.~Schuurmans, J.~D. Lafferty, C.~K.~I. Williams, and
  A.~Culotta, editors, {\em Advances in Neural Information Processing Systems
  22}, pages 781--789. Curran Associates, Inc., 2009.

\bibitem{Ji09}
Shuiwang Ji and Jieping Ye.
\newblock An accelerated gradient method for trace norm minimization.
\newblock In {\em Proceedings of the 26th Annual International Conference on
  Machine Learning}, ICML '09, pages 457--464, 2009.

\bibitem{Jojic10}
Vladimir Jojic, Stephen Gould, and Daphne Koller.
\newblock Accelerated dual decomposition for map inference.
\newblock In {\em Proceedings of the 27th International Conference on
  International Conference on Machine Learning}, ICML'10, pages 503--510, 2010.

\bibitem{Kantorovich42}
Leonid Kantorovich.
\newblock On the transfer of masses.
\newblock {\em Dokl. Akad. Nauk. SSSR}, 37(7-8):227--220, 1942.

\bibitem{Kantorovich48}
Leonid Kantorovich.
\newblock On a problem of monge.
\newblock {\em Uspekhi Mat. Nauk.}, 3:225--226, 1948.

\bibitem{Kich-ARMA}
S.~Kichenassamy, Arun Kumar, Peter Olver, Allen Tannenbaum, and Anthony Yezzi.
\newblock Conformal curvature flows: From phase transistions to active vision.
\newblock {\em Archive for Rational Mechanics and Analysis}, 134:275--301,
  1996.

\bibitem{Krichene15}
Walid Krichene, Alexandre Bayen, and Peter~L Bartlett.
\newblock Accelerated mirror descent in continuous and discrete time.
\newblock In C.~Cortes, N.~D. Lawrence, D.~D. Lee, M.~Sugiyama, and R.~Garnett,
  editors, {\em Advances in Neural Information Processing Systems 28}, pages
  2845--2853. Curran Associates, Inc., 2015.

\bibitem{Li15}
Huan Li and Zhouchen Lin.
\newblock Accelerated proximal gradient methods for nonconvex programming.
\newblock In C.~Cortes, N.~D. Lawrence, D.~D. Lee, M.~Sugiyama, and R.~Garnett,
  editors, {\em Advances in Neural Information Processing Systems 28}, pages
  379--387. Curran Associates, Inc., 2015.

\bibitem{MennucciNEW}
A.~Mennucci, A.~Yezzi, and G.~Sundaramoorthi.
\newblock Properties of sobolev-type metrics in the space of curves.
\newblock {\em Interfaces and Free Boundaries}, 10:423--445, 2008.

\bibitem{Monge}
G.~Monge.
\newblock M\'{e}moire sur la th\'{e}orie des d\'{e}blais et des remblais.
\newblock In {\em De l'Imprimerie Royale}. 1781.

\bibitem{Mukherjee13}
Indraneel Mukherjee, Kevin Canini, Rafael Frongillo, and Yoram Singer.
\newblock Parallel boosting with momentum.
\newblock In H.~Blockeel, K.~Kersting, S.~Nijssen, and F.~Zelezny, editors,
  {\em Machine Learning and Knowledge Discovery in Databases}, pages 17--32.
  Springer, Berlin, 2013.

\bibitem{Nesterov83}
Yurii Nesterov.
\newblock A method of solving a convex programming problem with convergence
  rate o (1/k2).
\newblock In {\em Soviet Mathematics Doklady}, volume~27, pages 372--376, 1983.

\bibitem{Nesterov05}
Yurii Nesterov.
\newblock Smooth minimization of non-smooth functions.
\newblock {\em Math. Program.}, 103(1):127--152, 2005.

\bibitem{Nesterov08}
Yurii Nesterov.
\newblock Accelerating the cubic regularization of newton's method on convex
  problems.
\newblock {\em Math. Program.}, 112(1):159--181, 2008.

\bibitem{Nesterov13}
Yurii Nesterov.
\newblock Gradient methods for minimizing composite functions.
\newblock {\em Math. Program.}, 140(1):125--161, 2013.

\bibitem{Nesterov14}
Yurii Nesterov.
\newblock {\em Introductory Lectures on Convex Optimization: A Basic Course}.
\newblock Springer Publishing Company, Incorporated, 1 edition, 2014.

\bibitem{Nesterov06}
Yurii Nesterov and Boris~T. Polyak.
\newblock Cubic regularization of newton method and its global performance.
\newblock {\em Math. Program.}, 108(1):177--205, 2006.

\bibitem{Niethammer03}
Marc Niethamer and Allen Tannenbaum.
\newblock Dynamic level sets for visual tracking.
\newblock In {\em Proceedings of IEEE Conference on Decision and Control},
  pages 4883--4888, 2003.

\bibitem{Niethammer04}
Marc Niethamer and Allen Tannenbaum.
\newblock Dynamic geodesic snakes for visual tracking.
\newblock In {\em IEEE Conference on Computer Vision and Pattern Recognition},
  pages 4883--4888, 2004.

\bibitem{Niethammer06}
Marc Niethamer and Allen Tannenbaum.
\newblock Dynamic geodesic snakes for visual tracking.
\newblock {\em IEEE Transactions on Automatic Control}, 51(4):562--579, 2006.

\bibitem{Niethammer08}
Marc Niethamer, Patricio Vela, and Allen Tannenbaum.
\newblock Geometric observers for dynamically evolving curves.
\newblock {\em IEEE Trans. Pattern Analysis and Machine Intelligence},
  30(6):1093--1108, 2008.

\bibitem{ODonoghue15}
Brendan O'Donoghue and Emmanuel Cand\`{e}s.
\newblock Adaptive restart for accelerated gradient schemes.
\newblock {\em Foundations of Computational Mechanics}, 15(3):715--732, 2015.

\bibitem{OsherParagios-book}
Stanely Osher and Nikos Paragios.
\newblock {\em Geometric Level Set Methods in Imaging, Vision and Graphics}.
\newblock Springer, New York, 2003.

\bibitem{OsherSethian88}
Stanely Osher and James Sethian.
\newblock Fronts propagation with curvature dependent speed: Algorithms based
  on hamilton-jacobi formulations.
\newblock {\em Journal of Computational Physics}, 79:12--49, 1988.

\bibitem{Pock09}
Thomas Pock, Antonin Chambolle, Daniel Cremers, and Horst Bischof.
\newblock A convex relaxation approach for computing minimal partitions.
\newblock In {\em Computer Vision and Pattern Recognition}, pages 810--817.
  IEEE, 2009.

\bibitem{Sapiro-book}
Guillermo Sapiro.
\newblock {\em Geometric Partial Differential Equations and Image Analysis}.
\newblock Cambridge Press, Cambridge, England, 2000.

\bibitem{Sethian-book}
James Sethian.
\newblock {\em Level Set Methods: Evolving Interfaces in Geometry, Fluid
  Mechanics, Computer Vision, and Material Science}.
\newblock Cambridge University Press, 1996.

\bibitem{Su14}
Weijie Su, Stephen Boyd, and Emmanuel Cand\`{e}s.
\newblock A differential equation for modeling nesterov's accelerated gradient
  method: Theory and insights.
\newblock In {\em Advances in Neural Information Processing Systems}, pages
  2510--2518, 2014.

\bibitem{Sundaramoorthi06}
G.~Sundaramoorthi, J.~Jackson, A.~Yezzi, and A.~Mennucci.
\newblock Tracking with sobolev active contours.
\newblock In {\em Proc. Computer Vision and Pattern Recognition}, pages
  674--680, 2006.

\bibitem{Sundaramoorthi-SIAM}
G.~Sundaramoorthi, A.~Mennucci, S.~Soatto, and A.~Yezzi.
\newblock A new geometric metric in the space of curves, and applications to
  tracking deforming objects by prediction and filtering.
\newblock {\em SIAM J. Imaging Sciences}, 4:109--145, 2011.

\bibitem{Sundaramoorthi-CDC}
G.~Sundaramoorthi, A.~Mennucci, A.~Yezzi, and S.~Soatto.
\newblock Tracking deforming objects by filtering and prediction in the space
  of curves.
\newblock In {\em Proc. IEEE Conf. Decision and Control}, pages 2395--2401,
  2009.

\bibitem{Sundaramoorthi05}
G.~Sundaramoorthi, A.~Yezzi, and A.~Mennucci.
\newblock Sobolev active contours.
\newblock In {\em Workshop on Variational Geometric and Level Set Methods in
  Computer Vision}, pages 109--120, 2005.

\bibitem{Sundaramoorthi07}
Ganesh Sundaramoorthi, Anthony Yezzi, and Andrea Mennucci.
\newblock Sobolev active contours.
\newblock {\em Int. J. Computer Vision}, 7:345--366, 2007.

\bibitem{Troutman96}
J.~L. Troutman.
\newblock {\em Variational Calclus and Optimal Control}.
\newblock Springer-Verlag, New York, 1996.

\bibitem{Villani03}
C.~Villani.
\newblock Topics in optimal transportation.
\newblock In {\em Graduate Studies in Mathematics 58}. AMS, Providence RI,
  2003.

\bibitem{Wibisono16}
Andre Wibisono, Ashia~C Wilson, and Michael~I Jordan.
\newblock A variational perspective on accelerated methods in optimization.
\newblock {\em Proceedings of the National Academy of Sciences}, page
  201614734, 2016.

\bibitem{Yang15}
Y.~Yang and Ganesh Sundaramoorthi.
\newblock Shape tracking with occlusions via coarse-to-fine region based
  sobolev descent.
\newblock {\em Trans. Pattern Analysis and Machine Intelligence}, 2015.

\end{thebibliography}

\end{document}